\documentclass[12pt, a4paper]{article}
\usepackage[latin1]{inputenc}
\usepackage{vmargin}
\usepackage{amsfonts}
\usepackage{amssymb}
\usepackage{graphicx}
\usepackage{amsmath}
\usepackage{amsthm}
\usepackage{setspace}
\usepackage{multirow}
\usepackage{hyperref}
\usepackage{cases}
\usepackage[bottom]{footmisc}
\setmarginsrb{1.2cm}{1.6cm}{1.2cm}{1.6cm}{0cm}{0cm}{0cm}{0.8cm}

\begin{document}

\newtheorem{Proposition}{Proposition}
\newtheorem{Theorem}{Theorem}
\newtheorem{Corollary}{Corollary}
\newtheorem{Def}{Definition}
\newtheorem{Lem}{Lemma}

\title{Mean field games equations with quadratic Hamiltonian: a specific approach\thanks{The author wishes to acknowledge the helpful conversations with Yves Achdou (Université Paris-Diderot), Pierre Cardaliaguet (Université Paris-Dauphine), Jean-Michel Lasry (Université Paris-Dauphine), Antoine Lemenant (Université Paris-Diderot), Pierre-Louis Lions (Collège de France), Stéphane Menozzi (Université Paris-Diderot), Vincent Millot (Université Paris-Diderot).}}
\author{Olivier Guéant\addtocounter{footnote}{5} \thanks{UFR de Math\'ematiques, Laboratoire Jacques-Louis Lions, Universit\'e Paris-Diderot. 175, rue du Chevaleret, 75013 Paris, France. \texttt{olivier.gueant@ann.jussieu.fr}}}
\date{}

\maketitle

\abstract{Mean field games models describing the limit of a large class of stochastic differential games, as the number of players goes to $+\infty$, have been introduced by J.-M. Lasry and P.-L. Lions in \cite{MFG1,MFG2,MFG3}. We use a change of variables to transform the mean field games (MFG) equations into a system of simpler coupled partial differential equations, in the case of a quadratic Hamiltonian. This system is then used to exhibit a monotonic scheme to build solutions of the MFG equations. Effective numerical methods based on this constructive scheme are presented and numerical experiments are carried out.}

\section*{Introduction}

Mean field games equations have been introduced by J.-M. Lasry and P.-L. Lions \cite{MFG1,MFG2,MFG3} to describe the dynamic equilibrium of stochastic differential games involving a continuum of players.\\

In the time-dependent case, these equations write:
$$\mathrm{(HJB)}\qquad  \partial_t u + \frac{\sigma^2}{2} \Delta u +H(\nabla u) = -f(x,m)$$
$$\mathrm{(K)}\qquad \partial_t m + \nabla \cdot (mH'(\nabla u))= \frac{\sigma^2}{2} \Delta m $$
with prescribed initial condition $m(0,\cdot) = m_0(\cdot) \ge 0$ and terminal condition $u(T,\cdot) = u_T(\cdot)$, where $u$ and $m$ are scalar functions defined on $[0,T]\times\Omega$, $\Omega$ typically being $(0,1)^d$.\\

In this paper, we focus on the particular case of quadratic hamiltonian $H(p) =\frac{p^2}{2}$. In this special case, a change of variables have been introduced by O. Guéant, J.-M. Lasry and P.-L. Lions in \cite{ParisPrinceton} to write the mean field games equations as two coupled heat equations with similar source terms. If indeed we introduce $\phi = \exp\left(\frac{u}{\sigma^2}\right)$ and $\psi = m \exp\left(-\frac{u}{\sigma^2}\right)$ then the system reduces to:

\begin{eqnarray*}
\partial_t \phi + \frac{\sigma^2}{2} \Delta \phi &=& -\frac{1}{\sigma^2} f(x,\phi\psi) \phi\\
\partial_t \psi - \frac{\sigma^2}{2} \Delta \psi &=& \frac{1}{\sigma^2} f(x,\phi\psi) \psi
\end{eqnarray*}
with $\phi(T,\cdot) = \exp\left(\frac{u_T(\cdot)}{\sigma^2}\right)$ and $\psi(0,\cdot) = \frac{m_0(\cdot)}{\phi(0,\cdot)}$.\\

We use this system to exhibit a constructive scheme for solutions to the mean field games equations. This constructive scheme starts with $\psi^0 =0$ and builds recursively two sequences $(\phi^{n+\frac{1}{2}})_n$ and $(\psi^{n+1})_n$ using the following equations:

\begin{eqnarray*}
\partial_t \phi^{n+\frac 12} + \frac{\sigma^2}{2} \Delta \phi^{n+\frac 12} &=& -\frac{1}{\sigma^2} f(x,\phi^{n+\frac 12}\psi^n) \phi^{n+\frac 12}\\
\partial_t \psi^{n+1} - \frac{\sigma^2}{2} \Delta \psi^{n+1} &=& \frac{1}{\sigma^2} f(x,\phi^{n+\frac 12}\psi^{n+1}) \psi^{n+1}
\end{eqnarray*}
with $\phi^{n+\frac 12}(T,\cdot) = \exp\left(\frac{u_T(\cdot)}{\sigma^2}\right)$ and $\psi^{n+1}(0,\cdot) = \frac{m_0(\cdot)}{\phi^{n+\frac 12}(0,\cdot)}$.\\

Then, $\phi$ and $\psi$ are obtained as the monotonic limit of the two sequences $(\phi^{n+\frac{1}{2}})_n$ and $(\psi^n)_n$ under usual assumptions of $f$.\\

Moreover, this scheme provides a new numerical method to solve mean field games that is completely different from the methods already proposed in the literature. Finite different schemes on $u$ and $m$ have been proposed, along with methods based on an equivalent formulation of the mean field game partial differential equations in the form of an optimization problem (see \cite{achdou2010capuzzo}, \cite{achdou2010mean}, \cite{lachapelle2010computation}, \cite{gueant2009reference} for the different proposals).\\
Here, we build a discrete counterpart to the above constructive scheme and the resulting algorithm to approximate $(\phi,\psi)$ consists in a sequence of totally implicit finite difference schemes that can be tackled with Newton methods.\\
The main difference with the preceding literature arises from the monotonicity properties of the scheme and from the peculiarity that the equilibrium $m$ is approached by a sequence that does not verify, except at the limit, the mass conservation principle.\\

In section 1, we recall the change of variables and derive the associated system of coupled parabolic equations. Then, section 2 is devoted to the introduction of the functional framework and we prove the main monotonicity properties of the system. Section 3 presents the constructive scheme and proves that we can have two monotonic sequences converging towards $\phi$ and $\psi$. Section 4 uses the same ideas as in the preceding sections, but in a discrete setting, to provide numerical schemes to approximate $\phi$ and $\psi$ numerically. Stability and convergence of the scheme are proved. Finally, section 5 presents the numerical experiments carried out and discusses its properties.\\

\section{From mean field games equations to a forward-backward system of heat equations with source terms}

We consider the mean field games equations introduced in \cite{MFG1,MFG2,MFG3} in the case of a quadratic Hamiltonian. These partial differential equations, hereafter denoted (MFG) are considered on the domain $[0,T]\times\Omega$, $\Omega$ standing for $(0,1)^d$, and consist in the following equations:

$$\mathrm{(HJB)}\qquad  \partial_t u + \frac{\sigma^2}{2} \Delta u +\frac{1}{2} |\nabla u|^2 = -f(x,m)$$
$$\mathrm{(K)}\qquad \partial_t m + \nabla \cdot (m\nabla u)= \frac{\sigma^2}{2} \Delta m $$

with:
\begin{itemize}
  \item Boundary conditions: $\frac{\partial u}{\partial \vec{n}} = \frac{\partial m}{\partial \vec{n}} = 0$ on $(0,T)\times\partial\Omega$
  \item Terminal condition: $u(T,\cdot) = u_T(\cdot)$ a given payoff whose regularity is to be precised
  \item Initial condition: $m(0,\cdot) = m_0(\cdot) \ge 0$ a given positive function in $L^1(\Omega)$, typically a probability distribution function.
\end{itemize}

The change of variables introduced in \cite{ParisPrinceton} is recalled in the following proposition:

\begin{Proposition}
Let us consider a smooth solution $(\phi,\psi)$ of the following system $(\mathcal{S})$ with $\phi > 0$:
\begin{eqnarray*}
\partial_t \phi + \frac{\sigma^2}{2} \Delta \phi &=& -\frac{1}{\sigma^2} f(x,\phi\psi) \phi \qquad (E_\phi)\\
\partial_t \psi - \frac{\sigma^2}{2} \Delta \psi &=& \frac{1}{\sigma^2} f(x,\phi\psi) \psi \qquad (E_\psi)
\end{eqnarray*}
with:
\begin{itemize}
  \item Boundary conditions: $\frac{\partial \phi}{\partial \vec{n}} = \frac{\partial \psi}{\partial \vec{n}} = 0$ on $(0,T)\times\partial\Omega$
  \item Terminal condition: $\phi(T,\cdot) = \exp\left(\frac{u_T(\cdot)}{\sigma^2}\right)$.
  \item Initial condition: $\psi(0,\cdot) = \frac{m_0(\cdot)}{\phi(0,\cdot)}$
\end{itemize}
Then $(u,m) = (\sigma^2 \ln(\phi),\phi\psi)$ defines a solution of (MFG).
\end{Proposition}

\noindent\textbf{Proof:}\\

Let us start with (HJB).

$$\partial_t u = \sigma^2 \frac{\partial_t \phi}{\phi}, \qquad \nabla u =  \sigma^2 \frac{\nabla \phi}{\phi} \qquad \Delta u =  \sigma^2 \frac{\Delta \phi}{\phi} - \sigma^2 \frac{|\nabla \phi|^2}{\phi^2}$$

Hence
\begin{eqnarray*}
\partial_t u + \frac{\sigma^2}{2} \Delta u +\frac{1}{2} |\nabla u|^2 &=& \sigma^2 \left[\frac{\partial_t \phi}{\phi} + \frac{\sigma^2}{2} \frac{\Delta \phi}{\phi} \right]\\
&=&\frac{\sigma^2}{\phi} \left[-\frac{1}{\sigma^2} f(x,\phi\psi) \phi\right]\\
&=&-f(x,m)\\
\end{eqnarray*}

Now, for the equation $(K)$:

$$
\partial_t m = \partial_t \phi \psi + \phi \partial_t \psi \qquad
\nabla \cdot(\nabla u m) =  \sigma^2 \nabla \cdot (\nabla \phi \psi) = \sigma^2 \left[\Delta \phi \psi + \nabla \phi \cdot \nabla \psi\right]$$
$$\Delta m = \Delta \phi \psi + 2 \nabla \phi \cdot \nabla \psi + \phi \Delta \psi$$

Hence
\begin{eqnarray*}
\partial_t m + \nabla \cdot (\nabla u m) &=& \partial_t \phi \psi + \phi \partial_t \psi + \sigma^2 \left[\Delta \phi \psi + \nabla \phi \cdot \nabla \psi\right]\\
&=& \psi \left[\partial_t \phi +  \sigma^2 \Delta \phi\right] + \phi \partial_t \psi + \sigma^2 \nabla \phi \cdot \nabla \psi\\
&=& \psi \left[\frac{\sigma^2}{2} \Delta \phi - \frac{1}{\sigma^2} f(x,\phi\psi) \phi\right] + \phi\left[\frac{\sigma^2}{2} \Delta \psi +\frac{1}{\sigma^2} f(x,\phi\psi) \psi\right] + \sigma^2 \nabla \phi \cdot \nabla \psi\\
&=& \frac{\sigma^2}{2} \Delta \phi \psi + \sigma^2 \nabla \phi \cdot \nabla \psi + \frac{\sigma^2}{2} \phi \Delta \psi\\
&=& \frac{\sigma^2}{2} \Delta m
\end{eqnarray*}

This proves the result since the boundary conditions and the initial and terminal conditions are coherent.\qed\\

Now, we will focus our attention on the study of the above system of equations $(\mathcal{S})$ and use it to design a constructive scheme for the couple $(\phi,\psi)$ and thus for the couple $(u,m)$ under regularity assumptions.\\

\section{Properties of $(\mathcal{S})$}

To study the system $(\mathcal{S})$ we introduce several hypotheses on $f$: we suppose that it is a decreasing function of its second variable, continuous in that variable and uniformly bounded. Moreover, to simplify the exposition\footnote{In terms of the initial mean field game problem, the optimal control $\nabla u$ and the subsequent distribution $m$ are not changed if we subtract $\|f\|_{\infty}$ to $f$.} we suppose that $f \le 0$.\\
It's noteworthy that the monotony hypothesis is to be linked to the usual proof of uniqueness for the mean field games equations (see \cite{MFG3}).\\

Now let's introduce the functional framework we are working in.

Let us note $\mathcal{P} \subset C([0,T],L^2(\Omega))$ the natural set for parabolic equations:
$$ g \in \mathcal{P} \iff g \in L^2(0,T,H^1(\Omega))\quad \mathrm{and } \quad \partial_t g \in L^2(0,T,H^{-1}(\Omega))$$

and let's introduce $\mathcal{P}_\epsilon = \lbrace g \in \mathcal{P}, g \ge \epsilon \rbrace$.

\begin{Proposition}
Suppose that $u_T \in L^\infty(\Omega)$.\\
$\forall \psi \in \mathcal{P}_0$, there is a unique weak solution $\phi$ to the following equation $(E_\phi)$:
$$\partial_t \phi + \frac{\sigma^2}{2} \Delta \phi = -\frac{1}{\sigma^2} f(x,\phi\psi) \phi \qquad (E_\phi)$$
with $\frac{\partial \phi}{\partial \vec{n}} = 0$ on $(0,T)\times\partial\Omega$ and $\phi(T,\cdot) = \exp\left(\frac{u_T(\cdot)}{\sigma^2}\right)$.\\

Hence $\Phi : \psi \in \mathcal{P}_0 \mapsto \phi \in \mathcal{P}$ is well defined.\\

Moreover, $\forall \psi \in \mathcal{P}_0,  \phi = \Phi(\psi) \in \mathcal{P}_\epsilon$ for $\epsilon = \exp\left(-\frac{1}{\sigma^2}\left(\|u_T\|_{\infty} + \|f\|_{\infty} T\right)\right)$
\end{Proposition}

\noindent\textbf{Proof:}\\

Let us consider $\psi \in \mathcal{P}_0$.\\

\underline{Existence of a weak solution $\phi$}:\\

Let us introduce $F_\psi : \varphi \in L^2(0,T,L^2(\Omega)) \mapsto \phi$ weak solution of the following linear parabolic equation:

$$\partial_t \phi + \frac{\sigma^2}{2} \Delta \phi = -\frac{1}{\sigma^2} f(x,\varphi\psi) \phi$$
with $\frac{\partial \phi}{\partial \vec{n}} = 0$ on $(0,T)\times\partial\Omega$ and $\phi(T,\cdot) = \exp\left(\frac{u_T(\cdot)}{\sigma^2}\right)$.\\

By classical linear parabolic equations theory, $\phi$  is in $\mathcal{P} \subset L^2(0,T,L^2(\Omega))$.\\

Our goal is to use Schauder's fixed-point theorem on $F_\psi$.\\

\emph{Compactness:}\\

Usual energy estimates (see \cite{evans}) give that there exists a constant $C$ that only depends on $\|u_T\|_{\infty}$, $\sigma$ and $\|f\|_{\infty}$ such that:

$$\forall (\psi,\varphi) \in \mathcal{P}_0 \times L^2(0,T,L^2(\Omega)) ,\quad  \|F_\psi(\varphi)\|_{L^2(0,T,H^1(\Omega))} + \|\partial_t F_\psi(\varphi)\|_{L^2(0,T,H^{-1}(\Omega))} \le C$$

Hence $F_\psi$ maps the closed ball $B_{L^2(0,T,L^2(\Omega))}(0,C)$ to a compact subset of $B_{L^2(0,T,L^2(\Omega))}(0,C)$.\\

\emph{Continuity:}\\

Let us now prove that $F_\psi$ is a continuous function.\\

Let us consider a sequence $(\varphi_n)_n$ of $L^2(0,T,L^2(\Omega))$ with $\varphi_n \to_{n \to \infty} \varphi$ in the $L^2(0,T,L^2(\Omega))$ sense.\\
Let us write $\phi_n = F_\psi(\varphi_n)$. We know from the above compactness result that we can extract from $(\phi_n)_n$ a new sequence denoted $(\phi_{n'})_{n'}$ that converges in the $L^2(0,T,L^2(\Omega))$ sense toward a function $\phi$. To prove that $F_\psi$ is continuous, we then need to show that $\phi$ cannot be different from $F_\psi(\varphi)$.\\

Now, because of the energy estimates, we know that $\phi$ is  in $\mathcal{P}$ and that we can extract another subsequence (still denoted $(\phi_{n'})_{n'}$), such that:

\begin{itemize}
  \item $\phi_{n'} \to \phi$ in the $L^2(0,T,L^2(\Omega))$ sense.
  \item $\nabla\phi_{n'} \rightharpoonup \nabla\phi$ weakly in $L^2(0,T,L^2(\Omega))$.
  \item $\partial_t \phi_{n'} \rightharpoonup \partial_t \phi$ weakly in $L^2(0,T,H^{-1}(\Omega))$.
\end{itemize}
and
\begin{itemize}
  \item $\varphi_{n'} \to \varphi$ almost everywhere.
\end{itemize}

By definition we have that $\forall w \in L^2(0,T,H^1(\Omega))$:
$$\int_0^T \langle \partial_t\phi_{n'}(t,\cdot), w(t,\cdot) \rangle_{H^{-1}(\Omega),H^1(\Omega)} dt-\frac{\sigma^2}{2}  \int_0^T \int_\Omega \nabla \phi_{n'}(t,x) \cdot \nabla w(t,x) dx dt$$
$$=- \frac{1}{\sigma^2} \int_0^T \int_\Omega f(x,\varphi_{n'}(t,x)\psi(t,x)) \phi_{n'}(t,x) w(t,x) dx dt$$

By weak convergence, the left hand side of the above equality converges toward

$$\int_0^T \langle \partial_t\phi(t,\cdot), w(t,\cdot) \rangle_{H^{-1}(\Omega),H^1(\Omega)} dt-\frac{\sigma^2}{2}  \int_0^T \int_\Omega \nabla \phi(t,x) \cdot \nabla w(t,x) dx dt$$

Since $f$ is a bounded continuous function and using dominated convergence theorem, the right hand side converges toward

$$- \frac{1}{\sigma^2} \int_0^T \int_\Omega f(x,\varphi(t,x)\psi(t,x)) \phi(t,x) w(t,x) dx dt$$

Hence $\phi = F_\psi(\varphi)$.\\

\emph{Schauder's Theorem}:\\

By Schauder's theorem, we then know that there exists a fixed-point $\phi$ to $F_\psi$ and hence a weak solution to the nonlinear parabolic equation $(E_\phi)$.\\

\underline{Positiveness of $\phi$}:\\

Let us consider a solution $\phi$ as above. If $I(t) =  \frac{1}{2}\int_\Omega (\phi(t,x)_-)^2 dx$, then:

\begin{eqnarray*}
I'(t) &=& -\int_\Omega \partial_t \phi(t,x) \phi(t,x)_-  dx\\
& = &  - \int_\Omega \left( \nabla\phi(t,x)\cdot\nabla(\phi(t,x)_-) - \frac{1}{\sigma^2} f(x,\phi(t,x)\psi(t,x)) \phi(t,x)\phi(t,x)_- \right) dx\\
& = &   - \int_\Omega \left(-|\nabla\phi(t,x)|^2 1_{\phi(t,x)\le0} + \frac{1}{\sigma^2} f(x,\phi(t,x)\psi(t,x)) (\phi(t,x)_-)^2 \right) dx\\
\end{eqnarray*}
\begin{eqnarray*}
& = & \int_\Omega |\nabla\phi(t,x)|^2 1_{\phi(t,x)\le0} - \int_\Omega \frac{1}{\sigma^2}  f(x,\phi(t,x)\psi(t,x)) (\phi(t,x)_-)^2 dx \\
& \ge & 0
\end{eqnarray*}
Since $I(T) = 0$ and $I\ge 0$, we know that $I=0$. Hence, $\phi$ is positive.\\

\underline{Uniqueness}:\\

Let us consider two weak solutions $\phi_1$ and $\phi_2$ to the equation $(E_\phi)$.\\

Let us introduce $J(t) =  \frac{1}{2}\int_\Omega (\phi_2(t,x)-\phi_1(t,x))^2 dx$. We have:

\begin{eqnarray*}
J'(t) &=& \int_\Omega (\partial_t \phi_2(t,x) - \partial_t \phi_1(t,x)) (\phi_2(t,x)-\phi_1(t,x)) dx\\
& = &  -\int_\Omega \frac{1}{\sigma^2} \left(f(x,\phi_2(t,x)\psi(t,x)) \phi_2(t,x) - f(x,\phi_1(t,x)\psi(t,x)) \phi_1(t,x)\right)(\phi_2(t,x) - \phi_1(t,x)) dx\\
& & +\int_\Omega |\nabla\phi_2(t,x)-\nabla\phi_1(t,x)|^2 dx \\
\end{eqnarray*}
Because of our assumptions on $f$, the function $ \xi \in \mathbb{R}_+ \mapsto \frac{1}{\sigma^2}f(x,\psi \xi) \xi$ is a decreasing function.\\

Hence, since $\phi_1$ and $\phi_2$ are positive, $J'(t) \ge 0$. Since $J(T) = 0$ and $J\ge 0$, we know that $J=0$. Hence, $\phi_1 = \phi_2$.\\

\underline{Lower bound to $\phi$}:\\

We can get a lower bound to $\phi$ through a subsolution taken as the solution of the following ordinary differential equation:

$$ \underline{\phi}'(t) = \frac{1}{\sigma^2} \|f\|_\infty \underline{\phi}(t) \qquad \underline{\phi}(T) = \exp\left(-\frac{\|u_T\|_\infty}{\sigma^2}\right)$$

Let us indeed consider $K(t) =  \frac{1}{2}\int_\Omega ((\underline{\phi}(t)-\phi(t,x))_+)^2 dx$. We have:
\begin{eqnarray*}
K'(t) &=& \int_\Omega (\underline{\phi}'(t) - \partial_t \phi(t,x)) (\underline{\phi}(t)-\phi(t,x))_+ dx\\
& = &  \int_\Omega \left(\frac{1}{\sigma^2} \|f\|_\infty \underline{\phi}(t) (\underline{\phi}(t)-\phi(t,x))_+ + |\nabla\phi(t,x)|^2 1_{\underline{\phi}(t)-\phi(t,x)\ge 0} \right.\\
&+&\left. \frac{1}{\sigma^2} \left(f(x,\phi(t,x)\psi(t,x)) \phi(t,x) (\underline{\phi}(t)-\phi(t,x))_+ \right)\right) dx\\
& \ge &  \frac{1}{\sigma^2} \int_\Omega   \left(\|f\|_\infty \underline{\phi}(t) + f(x,\phi(t,x)\psi(t,x)) \phi(t,x)\right) (\underline{\phi}(t)-\phi(t,x))_+ dx\\
& \ge &  \frac{1}{\sigma^2} \int_\Omega   \left(\|f\|_\infty +f(x,\phi(t,x)\psi(t,x))\right) \phi(t,x) (\underline{\phi}(t)-\phi(t,x))_+ dx\\
& \ge & 0
\end{eqnarray*}

Since $K(T) = 0$ and $K\ge 0$, we know that $K=0$.\\
Hence, $\phi(t,x) \ge \underline{\phi}(t) = e^{-\frac{\|u_T\|_\infty}{\sigma^2}}\exp\left( - \frac{1}{\sigma^2} \|f\|_\infty (T-t) \right)\ge \epsilon$ and the result follows.\qed\\

Now, we turn to a monotonicity result regarding $\Phi$.

\begin{Proposition}
$$\forall \psi_1 \le \psi_2 \in \mathcal{P}_0, \Phi(\psi_1) \ge \Phi(\psi_2) $$
\end{Proposition}

\noindent\textbf{Proof:}\\

Let us introduce $\phi_1 = \Phi(\psi_1)$ and $\phi_2 = \Phi(\psi_2)$.\\

Let us introduce $I(t) =  \frac{1}{2}\int_\Omega ((\phi_2(t,x)-\phi_1(t,x))_+)^2 dx$. We have:

\begin{eqnarray*}
I'(t) &=& \int_\Omega (\partial_t \phi_2(t,x) - \partial_t \phi_1(t,x)) (\phi_2(t,x)-\phi_1(t,x))_+ dx\\
& = &  \int_\Omega \bigg( |\nabla\phi_2(t,x)-\nabla\phi_1(t,x)|^2 1_{\phi_2(t,x)-\phi_1(t,x)\ge0}\\
&- &\left. \frac{1}{\sigma^2} \left(f(x,\phi_2(t,x)\psi_2(t,x)) \phi_2(t,x) - f(x,\phi_1(t,x)\psi_1(t,x)) \phi_1(t,x)\right)(\phi_2(t,x) - \phi_1(t,x))_+\right) dx\\
&\ge &\frac{1}{\sigma^2} \int_\Omega \left(f(x,\phi_1(t,x)\psi_1(t,x)) \phi_1(t,x) - f(x,\phi_2(t,x)\psi_2(t,x)) \phi_2(t,x)\right)(\phi_2(t,x) - \phi_1(t,x))_+ dx\\
&\ge &\frac{1}{\sigma^2} \int_\Omega \left(f(x,\phi_1(t,x)\psi_2(t,x)) \phi_1(t,x) - f(x,\phi_2(t,x)\psi_2(t,x)) \phi_2(t,x)\right)(\phi_2(t,x) - \phi_1(t,x))_+ dx\\
& \ge & 0
\end{eqnarray*}

Hence, since $I(T) = 0$ and $I\ge 0$, we know that $I=0$. Consequently, $\phi_1 \ge \phi_2$.\qed\\

We now turn to the second equation $(E_\psi)$ of the system $(\mathcal{S})$.

\begin{Proposition}
Let us fix $\epsilon > 0$ and suppose that $m_0 \in L^2(\Omega)$.\\
$\forall \phi \in \mathcal{P}_\epsilon$, there is a unique weak solution $\psi$ to the following equation $(E_\psi)$:
$$\partial_t \psi - \frac{\sigma^2}{2} \Delta \psi = \frac{1}{\sigma^2} f(x,\phi\psi) \psi \qquad (E_\psi)$$
with $\frac{\partial \psi}{\partial \vec{n}} = 0$ on $(0,T)\times\partial\Omega$ and $\psi(0,\cdot) = \frac{m_0(\cdot)}{\phi(0,\cdot)}$.\\

Hence $\Psi : \phi \in \mathcal{P}_\epsilon \mapsto \psi \in \mathcal{P}$ is well defined.\\

Moreover, $\forall \phi \in \mathcal{P}_\epsilon,  \psi = \Psi(\phi) \in \mathcal{P}_0$.
\end{Proposition}

\noindent\textbf{Proof:}\\

The proof of existence and uniqueness of a weak solution $\psi \in \mathcal{P}$ is the same as in Proposition 2. The only thing to notice is that the initial condition $\psi(0,\cdot)$ is in $\L^2(\Omega)$ because $m_0 \in L^2(\Omega)$ and $\phi$ bounded from below by $\epsilon > 0$.\\

Now, to prove that $\psi \ge 0$, let's introduce $I(t) =  \frac{1}{2}\int_\Omega (\psi(t,x)_-)^2 dx$, then:

\begin{eqnarray*}
I'(t) &=& -\int_\Omega \partial_t \psi(t,x) \psi(t,x)_-  dx\\
& = &  - \int_\Omega \left( -\nabla\psi(t,x)\cdot\nabla(\psi(t,x)_-) + \frac{1}{\sigma^2} f(x,\phi(t,x)\psi(t,x)) \psi(t,x)\psi(t,x)_- \right) dx\\
& = &   - \int_\Omega \left(|\nabla\psi(t,x)|^2 1_{\psi(t,x)\le0} - \frac{1}{\sigma^2} f(x,\phi(t,x)\psi(t,x)) (\psi(t,x)_-)^2 \right) dx\\
& \le & 0
\end{eqnarray*}
Since $I(0) = 0$ and $I\ge 0$, we know that $I=0$. Hence, $\psi$ is positive.\qed\\

Now, we turn to a monotonicity result regarding $\Psi$.

\begin{Proposition}
$$\forall \phi_1 \le \phi_2 \in \mathcal{P}_\epsilon, \Psi(\phi_1) \ge \Psi(\phi_2) $$
\end{Proposition}

\noindent\textbf{Proof:}\\

Let us introduce $\psi_1 = \Psi(\phi_1)$ and $\psi_2 = \Psi(\phi_2)$.\\

Let us introduce $I(t) =  \frac{1}{2}\int_\Omega ((\psi_2(t,x)-\psi_1(t,x))_+)^2 dx$. We have:

\begin{eqnarray*}
I'(t) &=& \int_\Omega (\partial_t \psi_2(t,x) - \partial_t \psi_1(t,x)) (\psi_2(t,x)-\psi_1(t,x))_+ dx\\
& = & \frac{1}{\sigma^2} \left(f(x,\phi_2(t,x)\psi_2(t,x)) \psi_2(t,x) - f(x,\phi_1(t,x)\psi_1(t,x)) \psi_1(t,x)\right)(\psi_2(t,x) - \psi_1(t,x))_+ dx\\
& -& \int_\Omega |\nabla\psi_2(t,x)-\nabla\psi_1(t,x)|^2 1_{\psi_2(t,x)-\psi_1(t,x)\ge0} dx\\
&\le &\frac{1}{\sigma^2} \int_\Omega \left(f(x,\phi_2(t,x)\psi_2(t,x)) \psi_2(t,x) - f(x,\phi_1(t,x)\psi_1(t,x)) \psi_1(t,x)\right)(\psi_2(t,x) - \psi_1(t,x))_+ dx\\
&\le &\frac{1}{\sigma^2} \int_\Omega \left(f(x,\phi_1(t,x)\psi_2(t,x)) \psi_2(t,x) - f(x,\phi_1(t,x)\psi_1(t,x)) \psi_1(t,x)\right)(\psi_2(t,x) - \psi_1(t,x))_+ dx\\
& \le & 0
\end{eqnarray*}

Now, $I(0) = \frac{1}{2} \int_\Omega{ m_0(x)\left(\left(\frac{1}{\phi_2(0,x)}-\frac{1}{\phi_1(0,x)}\right)_+\right)^2 dx} =0$. Hence since $I\ge 0$, we know that $I=0$. Consequently, $\psi_1 \ge \psi_2$.\qed\\

We will use these properties to design a constructive scheme for the couple $(\phi,\psi)$.

\section{A constructive scheme to solve the system $(\mathcal{S})$}

The scheme we consider involves two sequences $(\phi^{n+\frac{1}{2}})_n$ and $(\psi^n)_n$ that are built using the following recursive equations:

$$\psi^0 = 0$$
$$\forall n \in \mathbb{N}, \phi^{n+\frac{1}{2}} = \Phi(\psi^n)$$
$$\forall n \in \mathbb{N}, \psi^{n+1} = \Psi(\phi^{n+\frac{1}{2}}) $$

\begin{Theorem}
Suppose that $u_T \in L^\infty(\Omega)$ and that $m_0 \in L^2(\Omega)$.\\

Then, the above scheme has the following properties:\\

\begin{itemize}
  \item $(\phi^{n+\frac{1}{2}})_n$ is a decreasing sequence of $\mathcal{P}_\epsilon$ where $\epsilon$ is as in Proposition 2.
  \item $(\psi^n)_n$ is an increasing sequence of $\mathcal{P}_0$, bounded from above by $\Psi(\epsilon)$
  \item $(\phi^{n+\frac{1}{2}},\psi^n)_n$ converges for almost every $(t,x) \in (0,T)\times\Omega$, and in $L^2(0,T,L^2(\Omega))$, towards a couple $(\phi,\psi)$.
  \item $(\phi,\psi) \in \mathcal{P}_\epsilon\times\mathcal{P}_0$ is a weak solution of $(\mathcal{S})$.
\end{itemize}
\end{Theorem}

\noindent\textbf{Proof:}\\

By immediate induction we get from Propositions 2 and 4 that the two sequences are well defined and in the appropriate spaces.\\

Now, as far as monotonicity is concerned we have that $\psi^1 = \Psi(\phi^{\frac 12}) \ge 0 = \psi^0$. Hence, if for a given $n \in \mathbb{N}$ we have $\psi^{n+1} \ge \psi^n$, then Proposition 3 gives: $$\phi^{n+\frac{3}{2}} = \Phi(\psi^{n+1}) \le \Phi(\psi^n) = \phi^{n+\frac{1}{2}}$$
Applying now the function $\Psi$ we get:$$\psi^{n+2} = \Psi(\phi^{n+\frac{3}{2}}) \ge \Psi(\phi^{n+\frac{1}{2}}) = \psi^{n+1}$$
By induction we then have that $(\phi^{n+\frac{1}{2}})_n$ is decreasing and $(\psi^n)_n$ is increasing.\\

Moreover, since $\phi^{n+\frac{1}{2}} \ge \epsilon$, we have that $\psi^{n+1} = \Psi(\phi^{n+\frac{1}{2}}) \le \Psi(\epsilon)$.\\

Now this monotonic behavior allows to define two limit functions $\phi$ and $\psi$ in $L^2(0,T,L^2(\Omega))$ and the convergence is almost everywhere and in $L^2(0,T,L^2(\Omega))$.\\

Now, we want to show that $(\phi,\psi)$ is a weak solution of $(\mathcal{S})$ and to this purpose we use the energy estimates of the parabolic equations.\\

We know that there exists a constant $C>0$, that only depends on $\|u_T\|_{\infty}$, $\sigma$ and $\|f\|_{\infty}$ such that:

$$\forall n \in \mathbb{N}, \qquad \|\phi^{n+\frac{1}{2}}\|_{L^2(0,T,H^1(\Omega))} + \|\partial_t \phi^{n+\frac{1}{2}}\|_{L^2(0,T,H^{-1}(\Omega))} \le C$$

Hence, $\phi \in \mathcal{P}$ and we can extract a subsequence $(\phi^{n'+\frac{1}{2}})_{n'}$ such that:

\begin{itemize}
  \item $\phi^{n'+\frac{1}{2}} \to \phi$ in the $L^2(0,T,L^2(\Omega))$ sense and almost everywhere.
  \item $\nabla\phi^{n'+\frac{1}{2}} \rightharpoonup \nabla\phi$ weakly in $L^2(0,T,L^2(\Omega))$.
  \item $\partial_t \phi^{n'+\frac{1}{2}} \rightharpoonup \partial_t \phi$ weakly in $L^2(0,T,H^{-1}(\Omega))$.
\end{itemize}

Now, for $w \in L^2(0,T,H^1(\Omega))$
$$\int_0^T \langle \partial_t\phi^{n'+\frac{1}{2}}(t,\cdot), w(t,\cdot) \rangle_{H^{-1}(\Omega),H^1(\Omega)} dt -\frac{\sigma^2}{2}  \int_0^T \int_\Omega \nabla \phi^{n'+\frac{1}{2}}(t,x) \cdot \nabla w(t,x) dx dt$$
$$=- \frac{1}{\sigma^2} \int_0^T \int_\Omega f(x,\phi^{n'+\frac{1}{2}}(t,x)\psi^{n'}(t,x)) \phi^{n'+\frac{1}{2}}(t,x) w(t,x) dx dt$$

Using the weak convergence stated above, the continuity hypothesis on $f$ and the dominated convergence theorem for the last term we get that $\forall w \in L^2(0,T,H^1(\Omega))$
$$\int_0^T \langle \partial_t\phi(t,\cdot), w(t,\cdot) \rangle_{H^{-1}(\Omega),H^1(\Omega)} dt-\frac{\sigma^2}{2}  \int_0^T \int_\Omega \nabla \phi(t,x) \cdot \nabla w(t,x) dx dt$$
$$=- \frac{1}{\sigma^2} \int_0^T \int_\Omega f(x,\phi(t,x)\psi(t,x)) \phi(t,x) w(t,x) dx dt$$

Hence, for almost every $t\in (0,T)$ and $\forall v \in H^1(\Omega)$,
$$\langle \partial_t\phi(t,\cdot), v \rangle_{H^{-1}(\Omega),H^1(\Omega)} -\frac{\sigma^2}{2}  \int_\Omega \nabla \phi(t,x) \cdot \nabla v(x) dx $$
$$=- \frac{1}{\sigma^2} \int_\Omega f(x,\phi(t,x)\psi(t,x)) \phi(t,x) v(x) dx$$
and the terminal condition is obviously satisfied.\\

Obviously, we also have $\phi \ge \epsilon$.\\

Now, turning to the second equation, the proof works the same. The only additional thing to notice is that $$\psi(0,\cdot) = \lim_{n\to\infty}\psi^{n+1}(0,\cdot) =  \lim_{n\to\infty} \frac{m_0}{\phi^{n+\frac{1}{2}}(0,\cdot)} = \frac{m_0}{\phi(0,\cdot)}$$ where the limits are in the $L^2(\Omega)$ sense.\\

Hence $(\phi,\psi)$ is indeed a weak solution of $(\mathcal{S})$.\qed\\

We exhibited a way to build a solution to the system $(\mathcal{S})$. However, it should be noticed that the two functions $\phi$ and $\psi$ have been introduced to transform the initial (MFG) system into two simpler partial differential equations. The functions we are interested in, as far as the mean field game is concerned, are rather $u$ and $m$.\\
To well understand the nature of the change of variables and of the constructive scheme, let us introduce the sequence $(m^{n+1})_n$ where $m^{n+1} = \phi^{n+\frac 12} \psi^{n+1}$. From Theorem $1$, we know that $(m^{n+1})_n$ converges almost everywhere and in $L^1$ toward the function $m = \phi \psi$ for which we have conservation of mass along the trajectory. However, this property is not true for $m^{n+1}$ as it is stated in the following proposition:

\begin{Proposition}
Let us consider $n \in \mathbb{N}$ and let's denote $M^{n+1}(t) = \int_\Omega m^{n+1}(t,x) dx$ the total mass of $m^{n+1}$ at date $t$.\\ Then, there may be a loss of mass along the trajectory in the sense that:
$$\frac{d\ }{dt}M^{n+1}(t) = \int_\Omega \psi^{n+1}(t,x) \phi^{n+\frac 12}(t,x) \left(f\left(x,\psi^{n+1}(t,x) \phi^{n+\frac 12}(t,x)\right)-f\left(x,\psi^{n}(t,x) \phi^{n+\frac 12}(t,x)\right)\right) \le 0$$
\end{Proposition}

\noindent \textbf{Proof:}\\

From the regularity obtained above, we can write:

\begin{eqnarray*}
  \frac{d\ }{dt}M^{n+1}(t) &=& \langle \partial_t\phi^{n+\frac{1}{2}}(t,\cdot), \psi^{n+1}(t,\cdot) \rangle_{H^{-1}(\Omega),H^1(\Omega)} + \langle \partial_t\psi^{n+1}(t,\cdot), \phi^{n+\frac{1}{2}}(t,\cdot) \rangle_{H^{-1}(\Omega),H^1(\Omega)} \\
   &=& \frac{\sigma^2}{2}\int_\Omega \nabla \phi^{n+\frac{1}{2}}(t,x) \cdot \nabla \psi^{n+1}(t,x) dx\\
   &-& \frac{1}{\sigma^2} \int_\Omega \phi^{n+\frac{1}{2}}(t,x) \psi^{n+1}(t,x) f\left(x,\phi^{n+\frac{1}{2}}(t,x) \psi^{n}(t,x)\right) dx \\
   &-& \frac{\sigma^2}{2}\int_\Omega \nabla \phi^{n+\frac{1}{2}}(t,x) \cdot \nabla \psi^{n+1}(t,x) dx\\
   &+& \frac{1}{\sigma^2} \int_\Omega \phi^{n+\frac{1}{2}}(t,x) \psi^{n+1}(t,x) f\left(x,\phi^{n+\frac{1}{2}}(t,x) \psi^{n+1}(t,x)\right) dx \\
   &=&\int_\Omega \psi^{n+1}(t,x) \phi^{n+\frac 12}(t,x) \left(f\left(x,\psi^{n+1}(t,x) \phi^{n+\frac 12}(t,x)\right)-f\left(x,\psi^{n}(t,x) \phi^{n+\frac 12}(t,x)\right)\right)\\
   &\le& 0
\end{eqnarray*}
\qed\\

This property shows that the constructive scheme is rather original since it basically consists in building a probability distribution function using a sequence of functions in $L^1$ that only has the right total mass asymptotically. For the same reason, although there will be no systematic loss of mass, the numerical scheme we present in the next section consists in approximating probability distribution functions without taking care to mass conservation.

\section{A numerical scheme}

The constructive scheme of Theorem 1 allows to design a monotonic numerical scheme to solve the system $(\mathcal{S})$ and in turn to approximate $u$ and $m$.\\

For the sake of simplicity, we present the scheme with $\Omega = (0,1)$ but the results would be similar in higher dimensions. We also assume that $u_T$ and $m_0$ are bounded functions.\\

Let us consider a uniform subdivision $(t_0, \ldots, t_I)$ of $[0,T]$ where $t_i = i\Delta t$ for $i \in \mathcal{I} = \lbrace 0, \ldots, I \rbrace$ and a uniform subdivision $(x_0, \ldots, x_J)$ of the interval $[0,1]$ where $x_j = j\Delta x$ for $j \in \mathcal{J} = \lbrace 0, \ldots, J \rbrace$ .

We consider a recursive sequence of finite difference schemes in which $\widehat{\psi}^n_{i,j}$ and $\widehat{\phi}^{n+\frac 12}_{i,j}$ stand for the approximations of $\psi^n$ and $\phi^{n+\frac 12}$ at point $(t_i,x_j)$. For convenience, we also define $\widehat{\psi}^n_{i,-1} = \widehat{\psi}^n_{i,0}$ and $\widehat{\psi}^n_{i,J+1} = \widehat{\psi}^n_{i,J}$ to take account of Neumann conditions at the border.\\

The numerical scheme is as follows:\\

$$\qquad \widehat{\psi}^0_{i,j} = 0, \quad (i,j) \in \mathcal{I}\times\mathcal{J}$$
and for $n\ge0$:
$$\frac{\widehat{\phi}^{n+\frac 12}_{i+1,j} - \widehat{\phi}^{n+\frac 12}_{i,j}}{\Delta t} + \frac{\sigma^2}{2} \frac{\widehat{\phi}^{n+\frac 12}_{i,j+1} - 2\widehat{\phi}^{n+\frac 12}_{i,j} + \widehat{\phi}^{n+\frac 12}_{i,j-1}}{(\Delta x)^2} = - \frac{1}{\sigma^2} f(x_j,\widehat{\phi}^{n+\frac 12}_{i,j}\widehat{\psi}^{n}_{i,j})\widehat{\phi}^{n+\frac 12}_{i,j}, \quad i\in \mathcal{I}-\lbrace I \rbrace, j\in \mathcal{J}$$
$$\widehat{\phi}^{n+\frac 12}_{I,j} = \exp\left(\frac{u_T(x_j)}{\sigma^2}\right),\quad j\in \mathcal{J}$$

$$\frac{\widehat{\psi}^{n+1}_{i+1,j} - \widehat{\psi}^{n+1}_{i,j}}{\Delta t} - \frac{\sigma^2}{2} \frac{\widehat{\psi}^{n+1}_{i+1,j+1} - 2\widehat{\psi}^{n+1}_{i+1,j} + \widehat{\psi}^{n+1}_{i+1,j-1}}{(\Delta x)^2} = \frac{1}{\sigma^2} f(x_j,\widehat{\phi}^{n+\frac 12}_{i+1,j}\widehat{\psi}^{n+1}_{i+1,j})\widehat{\psi}^{n+1}_{i+1,j}, \quad i\in \mathcal{I}-\lbrace I \rbrace, j\in \mathcal{J}$$
$$\widehat{\psi}^{n+1}_{0,j} = \frac{m_0(x_j)}{\widehat{\phi}^{n+\frac 12}_{0,j}},\quad j\in \mathcal{J}$$

To study this scheme, let's introduce $\mathcal{M} = M_{I+1,J+1}(\mathbb{R})$ and $$\mathcal{M}_\epsilon = \lbrace (m_{i,j})_{i \in \mathcal{I},j \in \mathcal{J}}  \in \mathcal{M}, \quad \forall (i,j) \in \mathcal{I}\times\mathcal{J},  m_{i,j} \ge \epsilon \rbrace$$

Let us start with the discrete counterpart of Proposition 2.\\

\begin{Proposition}
$\forall \widehat{\psi} \in \mathcal{M}_0$, there is a unique solution $\widehat{\phi} \in \mathcal{M}$ to the following equation:
$$\frac{\widehat{\phi}_{i+1,j} - \widehat{\phi}_{i,j}}{\Delta t} + \frac{\sigma^2}{2} \frac{\widehat{\phi}_{i,j+1} - 2\widehat{\phi}_{i,j} + \widehat{\phi}_{i,j-1}}{(\Delta x)^2} = - \frac{1}{\sigma^2} f(x_j,\widehat{\phi}_{i,j}\widehat{\psi}_{i,j})\widehat{\phi}_{i,j}, \quad i\in \mathcal{I}-\lbrace I \rbrace, j\in \mathcal{J}$$
with $\widehat{\phi}_{I,j} = \exp\left(\frac{u_T(x_j)}{\sigma^2}\right),\quad \forall j\in \mathcal{J}$ and the conventions $\widehat{\phi}_{i,-1} = \widehat{\phi}_{i,0}$, $\widehat{\phi}_{i,J+1} = \widehat{\phi}_{i,J}$.\\

Hence $\widehat{\Phi} : \widehat{\psi} \in \mathcal{M}_0 \mapsto \widehat{\phi} \in \mathcal{M}$ is well defined.\\

Moreover, $\forall \widehat{\psi} \in \mathcal{M}_0,  \widehat{\phi} = \widehat{\Phi}(\widehat{\psi}) \in \mathcal{M}_\epsilon$ for $\epsilon = \exp\left(-\frac{1}{\sigma^2}\left(\|u_T\|_{\infty} + \|f\|_{\infty} T\right)\right)$
\end{Proposition}

\noindent\textbf{Proof:}\\

Let us consider $\widehat{\psi} \in \mathcal{M}_0$.\\

\underline{Existence and positiveness of a solution $\widehat{\phi}$}:\\

Let us introduce $\widehat{F}_{\widehat{\psi}} : \widehat{\varphi} \in \mathcal{M} \mapsto \widehat{\phi}$ solution of the following equations:

$$\frac{\widehat{\phi}_{i+1,j} - \widehat{\phi}_{i,j}}{\Delta t} + \frac{\sigma^2}{2} \frac{\widehat{\phi}_{i,j+1} - 2\widehat{\phi}_{i,j} + \widehat{\phi}_{i,j-1}}{(\Delta x)^2} = - \frac{1}{\sigma^2} f(x_j,\widehat{\varphi}_{i,j}\widehat{\psi}_{i,j})\widehat{\phi}_{i,j}, \quad i\in \mathcal{I}-\lbrace I \rbrace, j\in \mathcal{J}$$
with $\widehat{\phi}_{I,j} = \exp\left(\frac{u_T(x_j)}{\sigma^2}\right),\forall j\in \mathcal{J}$ and the conventions $\widehat{\phi}_{i,-1} = \widehat{\phi}_{i,0}$, $\widehat{\phi}_{i,J+1} = \widehat{\phi}_{i,J}$.\\

Such a matrix $\widehat{\phi}$ exists by straightforward backward induction. We indeed know $\widehat{\phi}_{I,j}, \forall j\in \mathcal{J}$, and if we know $\widehat{\phi}_{i+1,j}, \forall j\in \mathcal{J}$ for some $i\in \mathcal{I}-\lbrace I \rbrace$, then we have:

$$A\left( \begin{array}{c}
\widehat{\phi}_{i,0} \\
\vdots\\
\widehat{\phi}_{i,J}\end{array} \right) =  \left( \begin{array}{c}
\widehat{\phi}_{i+1,0} \\
\vdots\\
\widehat{\phi}_{i+1,J}\end{array} \right)$$ where:
$$A = \begin{pmatrix}
1+\Delta t \beta_{i,0}  & -\frac{\sigma^2\Delta t}{2(\Delta x)^2} &  0 & \cdots & 0\\
-\frac{\sigma^2\Delta t}{2(\Delta x)^2} & 1+\Delta t \beta_{i,1} & -\frac{\sigma^2\Delta t}{2(\Delta x)^2}  &  \ddots & \vdots\\
0 & \ddots & \ddots & \ddots  & 0\\
\vdots & \ddots   & -\frac{\sigma^2\Delta t}{2(\Delta x)^2} & 1+ \Delta t \beta_{i,J-1} & -\frac{\sigma^2\Delta t}{2(\Delta x)^2}\\
0 & \cdots & 0 & -\frac{\sigma^2\Delta t}{2(\Delta x)^2} & 1+ \Delta t \beta_{i,J} \\
\end{pmatrix}$$
where $\beta_{i,j} = \frac{\sigma^2}{(\Delta x)^2} - \frac{1}{\sigma^2} f(x_j,\widehat{\varphi}_{i,j}\widehat{\psi}_{i,j})$, $\forall i\in \mathcal{I}-\lbrace I\rbrace, \forall j \in \mathcal{J}-\lbrace 0,J\rbrace$, $\beta_{i,0} = \frac{\sigma^2}{2(\Delta x)^2} - \frac{1}{\sigma^2} f(x_0,\widehat{\varphi}_{i,0}\widehat{\psi}_{i,0})$ and $\beta_{i,J} = \frac{\sigma^2}{2(\Delta x)^2} - \frac{1}{\sigma^2} f(x_J,\widehat{\varphi}_{i,J}\widehat{\psi}_{i,J})$, $\forall i\in \mathcal{I}-\lbrace I\rbrace$.\\

Since $A$ is an $M$-matrix, $\widehat{\phi}_{i,j}$ is well defined $\forall j\in \mathcal{J}$, and by immediate induction $\widehat{\phi}$ is well defined.\\

Now, $\widehat{F}_{\widehat{\psi}}$ is well-defined and continuous since $f$ is continuous.\\

To prove that there is a fixed point to $\widehat{F}_{\widehat{\psi}}$, we will use Brouwer's theorem. The only thing that needs to be proved is that $\exists C > 0, \forall \widehat{\varphi} \in \mathcal{M}, \max_{(i,j) \in \mathcal{I}\times\mathcal{J}} \left\|\widehat{F}_{\widehat{\psi}}(\widehat{\varphi})_{i,j}\right\| \le C$. More precisely, we are going to prove that if $\widehat{\phi} = \widehat{F}_{\widehat{\psi}}(\widehat{\varphi})$ then:
$$\forall(i,j) \in \mathcal{I}\times\mathcal{J}, \quad 0 \le \widehat{\phi}_{i,j} \le  \max_{j \in \mathcal{J}} \widehat{\phi}_{I,j} \le \left\|\exp\left(\frac{u_T}{\sigma^2}\right)\right\|_{\infty}$$

Let us suppose first that $\underline{\widehat{\phi}} = \min_{(i,j) \in \mathcal{I}\times\mathcal{J} } \widehat{\phi}_{i,j} < 0$.\\
Then, let's consider $(i,j) \in \left\lbrace (i',j') \in \mathcal{I}\times\mathcal{J} | \widehat{\phi}_{i,j} = \underline{\widehat{\phi}}, \forall (i'',j'')\in \mathcal{I}\times\mathcal{J}, \widehat{\phi}_{i'',j''} = \underline{\widehat{\phi}} \Rightarrow i'' \le i' \right\rbrace$. We know that $i \neq I$ and hence:

$$0<\frac{\widehat{\phi}_{i+1,j} - \widehat{\phi}_{i,j}}{\Delta t} + \frac{\sigma^2}{2} \frac{\widehat{\phi}_{i,j+1} - 2\widehat{\phi}_{i,j} + \widehat{\phi}_{i,j-1}}{(\Delta x)^2} = - \frac{1}{\sigma^2} f(x_j,\widehat{\varphi}_{i,j}\widehat{\psi}_{i,j})\widehat{\phi}_{i,j}$$
Since $f\le 0$, this gives $\widehat{\phi}_{i,j} = \underline{\widehat{\phi}} \ge 0$ contrary to our hypothesis.\\
Hence, $\forall (i,j)\in \mathcal{I}\times\mathcal{J},  \widehat{\phi}_{i,j} \ge 0$.\\

Similarly let's suppose that $\overline{\widehat{\phi}} = \max_{(i,j)\in \mathcal{I}\times\mathcal{J}} \widehat{\phi}_{i,j} > \max_{j\in \mathcal{J}} \widehat{\phi}_{I,j}$. Then, let's consider a couple $(i,j) \in \left\lbrace (i',j') \in \mathcal{I}\times\mathcal{J} | \widehat{\phi}_{i,j} = \overline{\widehat{\phi}}, \forall (i'',j'') \in \mathcal{I}\times\mathcal{J}, \widehat{\phi}_{i'',j''} = \overline{\widehat{\phi}} \Rightarrow i'' \le i' \right\rbrace$. We have that $i \neq I$ and hence:

$$0 > \frac{\widehat{\phi}_{i+1,j} - \widehat{\phi}_{i,j}}{\Delta t} + \frac{\sigma^2}{2} \frac{\widehat{\phi}_{i,j+1} - 2\widehat{\phi}_{i,j} + \widehat{\phi}_{i,j-1}}{(\Delta x)^2} = - \frac{1}{\sigma^2} f(x_j,\widehat{\varphi}_{i,j}\widehat{\psi}_{i,j})\widehat{\phi}_{i,j}$$

But then, $\widehat{\phi}_{i,j} < 0$ and this is not possible. Hence, we indeed have that the maximum is attained for $i=I$, namely $\max_{(i,j)\in \mathcal{I}\times\mathcal{J}} \widehat{\phi}_{i,j} = \max_{j \in \mathcal{J}} \widehat{\phi}_{I,j}$.\\

Now, from Brouwer's theorem, we know that there exists $\widehat{\phi} \in \mathcal{M}_0$, such that $\widehat{F}_{\widehat{\psi}}(\widehat{\phi}) = \widehat{\phi}$.

\underline{Uniqueness of $\widehat{\phi}$}:\\

Let us suppose that there are two fixed points $\widehat{\phi}$ and $\widehat{\phi}'$ to $\widehat{F}_{\widehat{\psi}}$ and imagine for instance that $\exists (i,j) \in \mathcal{I}\times\mathcal{J}$ such that $\widehat{\phi}_{i,j}>\widehat{\phi}'_{i,j}$. Then, let's consider $\widehat{\delta} = \widehat{\phi} - \widehat{\phi}'$, $\overline{\widehat{\delta}} = \max_{(i,j)\in \mathcal{I}\times\mathcal{J}} \widehat{\delta}_{i,j} $ and $(i,j) \in \left\lbrace (i',j') \in \mathcal{I}\times\mathcal{J} | \widehat{\delta}_{i,j} = \overline{\widehat{\delta}}, \forall (i'',j'') \in \mathcal{I}\times\mathcal{J}, \widehat{\phi}_{i'',j''} = \overline{\widehat{\delta}} \Rightarrow i'' \le i' \right\rbrace$. We know that $i \neq I$ and hence:

$$0>\frac{\widehat{\delta}_{i+1,j} - \widehat{\delta}_{i,j}}{\Delta t} + \frac{\sigma^2}{2} \frac{\widehat{\delta}_{i,j+1} - 2\widehat{\delta}_{i,j} + \widehat{\delta}_{i,j-1}}{(\Delta x)^2} $$$$= - \frac{1}{\sigma^2} \left[f(x_j,\widehat{\phi}_{i,j}\widehat{\psi}_{i,j})\widehat{\phi}_{i,j} - f(x_j,\widehat{\phi'}_{i,j}\widehat{\psi}_{i,j})\widehat{\phi'}_{i,j} \right] $$
But since $\widehat{\phi}_{i,j} > \widehat{\phi}'_{i,j} \ge 0$ as above, the right hand side of the above equation must be greater than $0$, in contradiction with the above inequality.\\
Hence, $\widehat{\phi} \le \widehat{\phi}'$ and, by symmetry, $\widehat{\phi} = \widehat{\phi}'$. As a consequence  $\widehat{\Phi} : \widehat{\psi} \in \mathcal{M}_0 \mapsto \widehat{\phi} \in \mathcal{M}$ is well defined.\\

\underline{Lower bound for $\widehat{\phi}$}:\\

Let us consider $\underline{\widehat{\phi}}_{i,j} = \exp\left(-\frac{\|u_T\|_\infty}{\sigma^2}\right) \frac{1}{\left(1+ \Delta t \frac{\|f\|_{\infty}}{\sigma^2}\right)^{I-i}}$ solution of the following problem:

$$\frac{\underline{\widehat{\phi}}_{i+1,j} - \underline{\widehat{\phi}}_{i,j}}{\Delta t} + \frac{\sigma^2}{2} \frac{\underline{\widehat{\phi}}_{i,j+1} - 2\underline{\widehat{\phi}}_{i,j} + \underline{\widehat{\phi}}_{i,j-1}}{(\Delta x)^2} = \frac{1}{\sigma^2} \|f\|_{\infty} \underline{\widehat{\phi}}_{i,j}, \quad i\in \mathcal{I}-\lbrace I \rbrace, j\in \mathcal{J}$$
with $\underline{\widehat{\phi}}_{I,j} = \exp\left(\frac{-\|u_T\|_{\infty}}{\sigma^2}\right),\forall j\in \mathcal{J}$ and the conventions $\underline{\widehat{\phi}}_{i,-1} = \underline{\widehat{\phi}}_{i,0}$, $\underline{\widehat{\phi}}_{i,J+1} = \underline{\widehat{\phi}}_{i,J}$.\\

Now, let's write $\widehat{\delta} = \underline{\widehat{\phi}} - \widehat{\phi}$ and let's suppose that $\overline{\widehat{\delta}} = \max_{(i,j)\in \mathcal{I}\times\mathcal{J}} \widehat{\delta}_{i,j} > 0$. Then, let's consider $(i,j) \in \left\lbrace (i',j') \in \mathcal{I}\times\mathcal{J} | \widehat{\delta}_{i,j} = \overline{\widehat{\delta}}, \forall (i'',j'') \in \mathcal{I}\times\mathcal{J}, \widehat{\delta}_{i'',j''} = \overline{\widehat{\delta}} \Rightarrow i'' \le i' \right\rbrace$. By construction, we know that $i \neq I$ and hence:

$$0 > \frac{\widehat{\delta}_{i+1,j} - \widehat{\delta}_{i,j}}{\Delta t} + \frac{\sigma^2}{2} \frac{\widehat{\delta}_{i,j+1} - 2\widehat{\delta}_{i,j} + \widehat{\delta}_{i,j-1}}{(\Delta x)^2} = \frac{1}{\sigma^2}\left[ \|f\|_{\infty} \underline{\widehat{\phi}}_{i,j} + f(x_j,\widehat{\phi}_{i,j}\widehat{\psi}_{i,j})\widehat{\phi}_{i,j}\right]$$

But then, $\underline{\widehat{\phi}}_{i,j} < \widehat{\phi}_{i,j}$ in contradiction with the hypothesis.\\

In conclusion, we have $\forall (i,j) \in \mathcal{I}\times\mathcal{J}$:

$$\widehat{\phi}_{i,j} \ge \underline{\widehat{\phi}}_{i,j} \ge e{-\frac{\|u_T\|_\infty}{\sigma^2}} \frac{1}{\left(1+ \Delta t \frac{\|f\|_{\infty}}{\sigma^2}\right)^{I-i}} \ge e^{-\frac{\|u_T\|_\infty}{\sigma^2}} \frac{1}{\left(1+ \Delta t \frac{\|f\|_{\infty}}{\sigma^2}\right)^{\frac{T}{\Delta t}}} \ge \epsilon $$\qed\\

Let us now turn to a monotonicity property of the function $\widehat{\Phi}$.

\begin{Proposition}
$$\forall \widehat{\psi}_1 \le \widehat{\psi}_2 \in \mathcal{M}_0, \widehat{\Phi}(\widehat{\psi}_1) \ge \widehat{\Phi}(\widehat{\psi}_2)$$
\end{Proposition}

\noindent\textbf{Proof:}\\

Let us consider $\widehat{\psi}_1 \le \widehat{\psi}_2 \in \mathcal{M}_0$ and let's denote $\widehat{\phi}_1 = \widehat{\Phi}(\widehat{\psi}_1)$ and $\widehat{\phi}_2 = \widehat{\Phi}(\widehat{\psi}_2)$.\\

Now, let's write $\widehat{\delta} = \widehat{\phi}_2 - \widehat{\phi}_1$ and let's suppose that $\overline{\widehat{\delta}} = \max_{(i,j)\in \mathcal{I}\times\mathcal{J}} \widehat{\delta}_{i,j} > 0$. Then, let's consider $(i,j) \in \left\lbrace (i',j') \in \mathcal{I}\times\mathcal{J} | \widehat{\delta}_{i,j} = \overline{\widehat{\delta}}, \forall (i'',j'') \in \mathcal{I}\times\mathcal{J}, \widehat{\delta}_{i'',j''} = \overline{\widehat{\delta}} \Rightarrow i'' \le i' \right\rbrace$. We know that $i \neq I$ and hence:

$$0 > \frac{\widehat{\delta}_{i+1,j} - \widehat{\delta}_{i,j}}{\Delta t} + \frac{\sigma^2}{2} \frac{\widehat{\delta}_{i,j+1} - 2\widehat{\delta}_{i,j} + \widehat{\delta}_{i,j-1}}{(\Delta x)^2} = -\frac{1}{\sigma^2}\left[ f(x_j,\widehat{\phi}_{2,i,j}\widehat{\psi}_{2,i,j})\widehat{\phi}_{2,i,j} - f(x_j,\widehat{\phi}_{1,i,j}\widehat{\psi}_{1,i,j})\widehat{\phi}_{1,i,j}\right]$$

But then, $f(x_j,\widehat{\phi}_{2,i,j}\widehat{\psi}_{2,i,j})\widehat{\phi}_{2,i,j} - f(x_j,\widehat{\phi}_{1,i,j}\widehat{\psi}_{1,i,j})\widehat{\phi}_{1,i,j} > 0$ in contradiction with the hypothesis for the chosen $(i,j)$.\\

Hence, $\widehat{\phi}_1 \ge \widehat{\phi}_2$.\qed\\

Let us now turn to the second equation. We have exactly the same result as above with the same proof:

\begin{Proposition}
$\forall \widehat{\phi} \in \mathcal{M}_\epsilon$, there is a unique solution $\widehat{\psi} \in \mathcal{M}$ to the following equation:
$$\frac{\widehat{\psi}_{i+1,j} - \widehat{\psi}_{i,j}}{\Delta t} - \frac{\sigma^2}{2} \frac{\widehat{\psi}_{i+1,j+1} - 2\widehat{\psi}_{i+1,j} + \widehat{\psi}_{i+1,j-1}}{(\Delta x)^2} = \frac{1}{\sigma^2} f(x_j,\widehat{\phi}_{i+1,j}\widehat{\psi}_{i+1,j})\widehat{\psi}_{i+1,j}, \quad 0\le i \le I-1, 0\le j \le J$$
with $\widehat{\psi}_{0,j} = \frac{m_0(x_j)}{\widehat{\phi}_{0,j}},\forall j \in \mathcal{J}$ and the conventions $\widehat{\psi}_{i,-1} = \widehat{\psi}_{i,0}$, $\widehat{\psi}_{i,J+1} = \widehat{\psi}_{i,J}$.\\

Hence $\widehat{\Psi} : \widehat{\phi} \in \mathcal{M}_\epsilon \mapsto \widehat{\psi} \in \mathcal{M}$ is well defined.\\

Moreover, $\forall \widehat{\phi} \in \mathcal{M}_\epsilon,  \widehat{\psi} = \widehat{\Psi}(\widehat{\phi}) \in \mathcal{M}_0$.
\end{Proposition}

Similarly, the same monotonicity result can be proved:

\begin{Proposition}
$$\forall \widehat{\phi}_1 \le \widehat{\phi}_2 \in \mathcal{M}_\epsilon, \widehat{\Psi}(\widehat{\phi}_1) \ge \widehat{\Psi}(\widehat{\phi}_2)$$
\end{Proposition}

Now, using the functions introduced in the above propositions, we can write the numerical scheme in a more compact fashion:

$$\widehat{\psi}^0 = 0$$
$$\forall n \in \mathbb{N}, \widehat{\phi}^{n+\frac{1}{2}} = \widehat{\Phi}\left(\widehat{\psi}^n\right)$$
$$\forall n \in \mathbb{N}, \widehat{\psi}^{n+1} = \widehat{\Psi}\left(\widehat{\phi}^{n+\frac{1}{2}}\right) $$

Using the monotonicity properties, we get that the sequences $\left(\widehat{\psi}^{n}\right)_n$ and $\left(\widehat{\phi}^{n+\frac 12}\right)_n$ converge monotonically. More precisely we have that:

\begin{Proposition}
The numerical scheme verifies the following properties:
\begin{itemize}
  \item $\left(\widehat{\phi}^{n+\frac{1}{2}}\right)_n$ is a decreasing sequence of $\mathcal{M}_\epsilon$ where $\epsilon$ is as in Proposition 7.
  \item $\left(\widehat{\psi}^n\right)_n$ is an increasing sequence of $\mathcal{M}_0$, bounded from above by $\frac{\|m_0\|_{\infty}}{\epsilon}$.
  \item $\left(\widehat{\phi}^{n+\frac{1}{2}},\widehat{\psi}^n\right)_n$ converges towards a couple $(\widehat{\phi},\widehat{\psi})  \in \mathcal{M}_\epsilon\times\mathcal{M}_0$ that verifies:
$$\frac{\widehat{\phi}_{i+1,j} - \widehat{\phi}_{i,j}}{\Delta t} + \frac{\sigma^2}{2} \frac{\widehat{\phi}_{i,j+1} - 2\widehat{\phi}_{i,j} + \widehat{\phi}_{i,j-1}}{(\Delta x)^2} = - \frac{1}{\sigma^2} f(x_j,\widehat{\phi}_{i,j}\widehat{\psi}_{i,j})\widehat{\phi}_{i,j}, \quad i\in \mathcal{I}-\lbrace I \rbrace, j\in \mathcal{J}$$
$$\widehat{\phi}_{I,j} = \exp\left(\frac{u_T(x_j)}{\sigma^2}\right),\forall j\in \mathcal{J}$$

$$\frac{\widehat{\psi}_{i+1,j} - \widehat{\psi}_{i,j}}{\Delta t} - \frac{\sigma^2}{2} \frac{\widehat{\psi}_{i+1,j+1} - 2\widehat{\psi}_{i+1,j} + \widehat{\psi}_{i+1,j-1}}{(\Delta x)^2} = \frac{1}{\sigma^2} f(x_j,\widehat{\phi}_{i+1,j}\widehat{\psi}_{i+1,j})\widehat{\psi}_{i+1,j}, \quad i\in \mathcal{I}-\lbrace I \rbrace, j\in \mathcal{J}$$
$$\widehat{\psi}_{0,j} = \frac{m_0(x_j)}{\widehat{\phi}_{0,j}},\forall j\in \mathcal{J}$$
\end{itemize}
\end{Proposition}

\noindent\textbf{Proof:}\\

Let us prove by induction that $$\forall n \in \mathbb{N}, 0 \le \widehat{\psi}^{n} \le \widehat{\psi}^{n+1} \le \widehat{\Psi}(\epsilon), \qquad \epsilon \le \widehat{\phi}^{n+\frac{3}{2}} \le \widehat{\phi}^{n+\frac{1}{2}}$$
First, for $n=0$ we have that: $$\widehat{\psi}^1 = \widehat{\Psi}(\widehat{\Phi}(\widehat{\psi}^0)) \ge 0 = \widehat{\psi}^0 $$
Consequently, $$\widehat{\phi}^{\frac 12} = \widehat{\Phi}(\widehat{\psi}^0) \ge \widehat{\Phi}(\widehat{\psi}^1) = \widehat{\phi}^{\frac 32}$$
and because of Proposition 7, we know that $\widehat{\phi}^{\frac 12} \ge \widehat{\phi}^{\frac 32} \ge \epsilon$. Hence, $\widehat{\psi}^1 = \widehat{\Psi}(\widehat{\psi}^{\frac 12}) \le \widehat{\Psi}(\epsilon)$.\\

Now, let's suppose that the result is true for some $n \in \mathbb{N}$, then:

$$\widehat{\psi}^{n} \le \widehat{\psi}^{n+1} \Rightarrow \widehat{\phi}^{n+\frac 12} = \widehat{\Phi}(\widehat{\psi}^{n}) \ge \widehat{\Phi}(\widehat{\psi}^{n+1}) = \widehat{\phi}^{n+\frac 32}$$
$$\Rightarrow \widehat{\psi}^{n+1} = \widehat{\Psi}(\widehat{\phi}^{n+ \frac 12}) \le \widehat{\Psi}(\widehat{\phi}^{n+\frac 32}) = \widehat{\psi}^{n+2}$$
Now, by Proposition 7, we know that $\widehat{\phi}^{n+\frac 32} \ge \epsilon$ and thus $\widehat{\psi}^{n+2} = \widehat{\Psi}(\widehat{\psi}^{n+\frac 32}) \le \widehat{\Psi}(\epsilon)$ and with the same discrete maximum principle as above, we have that $\widehat{\Psi}(\epsilon)$ is bounded by $\frac{\|m_0\|_{\infty}}{\epsilon}$.

Consequently, $\forall (i,j) \in \mathcal{I}\times\mathcal{J}$, $(\widehat{\phi}_{i,j}^{n+\frac{1}{2}})_n$ is a decreasing sequence bounded from below by $\epsilon$ that thus converges towards a value denoted $\widehat{\phi}_{i,j}$. Similarly, $\forall (i,j) \in \mathcal{I}\times\mathcal{J}$, $(\widehat{\psi}_{i,j}^{n})_n$ is an increasing sequence bounded from above that thus converges towards a value denoted $\widehat{\psi}_{i,j}$.\\
The resulting couple $(\widehat{\phi},\widehat{\psi}) \in \mathcal{M}_\epsilon\times\mathcal{M}_0$ straightforwardly verifies the equations of the proposition.\qed\\

Now, to study convergence, let's introduce a norm $\|\cdot\|$ on the sets $\mathcal{M}$ by:

$$\forall m=(m_{i,j})_{i,j} \in \mathcal{M}, \|m\|^2 = \sup_{0\le i \le I} \frac{1}{J+1}\sum_{j=0}^J m^2_{i,j}$$

Let us also define the consistency errors $\tilde{\eta}^{n+\frac{1}{2}}$ and $\tilde{\eta}^{n+1}$ associated to the equations that define respectively $\phi^{n+\frac 12}$ and $\psi^{n+1}$:

$$\tilde{\eta}^{n+\frac{1}{2}}_{i,j}=\frac{\tilde{\phi}^{n+\frac 12}_{i+1,j} - \tilde{\phi}^{n+\frac 12}_{i,j}}{\Delta t} + \frac{\sigma^2}{2} \frac{\tilde{\phi}^{n+\frac 12}_{i,j+1} - 2\tilde{\phi}^{n+\frac 12}_{i,j} + \tilde{\phi}^{n+\frac 12}_{i,j-1}}{(\Delta x)^2} + \frac{1}{\sigma^2} f(x_j,\tilde{\phi}^{n+\frac 12}_{i,j}\tilde{\psi}^{n}_{i,j})\tilde{\phi}^{n+\frac 12}_{i,j}$$
$$\tilde{\eta}^{n+1}_{i,j}=\frac{\tilde{\psi}^{n+1}_{i+1,j} - \tilde{\psi}^{n+1}_{i,j}}{\Delta t} - \frac{\sigma^2}{2} \frac{\tilde{\psi}^{n+1}_{i+1,j+1} - 2\tilde{\psi}^{n+1}_{i+1,j} + \tilde{\psi}^{n+1}_{i+1,j-1}}{(\Delta x)^2} - \frac{1}{\sigma^2} f(x_j,\tilde{\phi}^{n+\frac 12}_{i+1,j}\tilde{\psi}^{n+1}_{i+1,j})\tilde{\psi}^{n+1}_{i+1,j}$$
where $\tilde{\phi}^{n+\frac 12}_{i,j} = \phi^{n+\frac 12}(t_i,x_j)$ and $\tilde{\psi}^{n+1}_{i,j} = \psi^{n+1}(t_i,x_j)$.\\

We are now ready to enounce a theorem that gives recursive stability bounds for the scheme.\\

\begin{Theorem}
Let us suppose, in addition to the hypotheses made above, that $f$ is uniformly Lipschitz with respect to its second variable, i.e.
$$\exists K, \forall x \in (0,1), \forall \xi_1, \xi_2 \in \mathbb{R}_+, |f(x,\xi_2) - f(x,\xi_1)| \le K |\xi_2 - \xi_1|$$
Then, $\forall \nu > 0$, $\exists C >0$, $\forall I,J \in \mathbb{N}$ such that $\frac 1{\Delta t}> 1+\frac{K}{\sigma^2}\max\left(\left\|e^{\frac{u_T}{\sigma^2}} \right\|^2_\infty,\frac{\|m_0\|^2_\infty}{\epsilon^2}\right)+\nu$, we have $\forall n \in \mathbb{N}$:
$$\|\widehat{\phi}^{n+\frac 12} - \tilde{\phi}^{n+\frac 12}\|^2 \le C \|\widehat{\psi}^{n} - \tilde{\psi}^{n}\|^2 + C \|\tilde{\eta}^{n+\frac{1}{2}}\|^2$$
$$\|\widehat{\psi}^{n+1} - \tilde{\psi}^{n+1}\|^2 \le C \|\widehat{\phi}^{n+\frac 12} - \tilde{\phi}^{n+\frac 12}\|^2 + C \|\tilde{\eta}^{n+1}\|^2$$
\end{Theorem}

\noindent\textbf{Proof:}\\

Let us denote for $n\in \mathbb{N}$, $\delta^{n + \frac 12}_{i,j} = \widehat{\phi}^{n + \frac 12}_{i,j} -  \tilde{\phi}^{n + \frac 12}_{i,j}$ and $\delta^{n}_{i,j} = \widehat{\psi}^{n}_{i,j} -  \tilde{\psi}^{n}_{i,j}$.\\

By definition we have for $i\in \mathcal{I}-\lbrace I \rbrace, j\in \mathcal{J}$:

$$\frac{\delta^{n+\frac 12}_{i+1,j} - \delta^{n+\frac 12}_{i,j}}{\Delta t} + \frac{\sigma^2}{2} \frac{\delta^{n+\frac 12}_{i,j+1} - 2\delta^{n+\frac 12}_{i,j} + \delta^{n+\frac 12}_{i,j-1}}{(\Delta x)^2} = - \frac{1}{\sigma^2} \left(f(x_j,\widehat{\phi}^{n+\frac 12}_{i,j}\widehat{\psi}^n_{i,j})\widehat{\phi}^{n+\frac 12}_{i,j} - f(x_j,\tilde{\phi}^{n+\frac 12}_{i,j}\tilde{\psi}^n_{i,j})\tilde{\phi}^{n+\frac 12}_{i,j}\right) - \tilde{\eta}^{n + \frac 12}_{i,j}$$

Hence:

$$\left(\delta^{n+\frac 12}_{i,j}\right)^2 = \delta^{n+\frac 12}_{i+1,j} \delta^{n+\frac 12}_{i,j} + \frac{\sigma^2 \Delta t}{2 (\Delta x)^2} \delta^{n+\frac 12}_{i,j} \left(\delta^{n+\frac 12}_{i,j+1} - 2\delta^{n+\frac 12}_{i,j} + \delta^{n+\frac 12}_{i,j-1}\right)$$
$$+ \frac{\Delta t}{\sigma^2} \delta^{n+\frac 12}_{i,j} \left(f(x_j,\widehat{\phi}^{n+\frac 12}_{i,j}\widehat{\psi}^n_{i,j})\widehat{\phi}^{n+\frac 12}_{i,j} - f(x_j,\tilde{\phi}^{n+\frac 12}_{i,j}\tilde{\psi}^n_{i,j})\tilde{\phi}^{n+\frac 12}_{i,j}\right) + \delta^{n+\frac 12}_{i,j} \tilde{\eta}^{n+\frac 12}_{i,j} \Delta t$$
$$\Rightarrow \frac 12 \left(\delta^{n+\frac 12}_{i,j}\right)^2 \le \frac 12 \left(\delta^{n+\frac 12}_{i+1,j}\right)^2 + \frac{\sigma^2 \Delta t}{2 (\Delta x)^2} \delta^{n+\frac 12}_{i,j} \left(\delta^{n+\frac 12}_{i,j+1} - 2\delta^{n+\frac 12}_{i,j} + \delta^{n+\frac 12}_{i,j-1}\right)$$
$$+ \frac{\Delta t}{\sigma^2} \delta^{n+\frac 12}_{i,j} \left(f(x_j,\widehat{\phi}^{n+\frac 12}_{i,j}\widehat{\psi}^n_{i,j})\widehat{\phi}^{n+\frac 12}_{i,j} - f(x_j,\tilde{\phi}^{n+\frac 12}_{i,j}\tilde{\psi}^n_{i,j})\tilde{\phi}^{n+\frac 12}_{i,j}\right) + \frac {\Delta t}{2} \left(\delta^{n+\frac 12}_{i,j}\right)^2 + \frac {\Delta t}{2} \left(\tilde{\eta}^{n+\frac 12}_{i,j}\right)^2$$

Now,

\begin{eqnarray*}
   & & \delta^{n+\frac 12}_{i,j} \left(f(x_j,\widehat{\phi}^{n+\frac 12}_{i,j}\widehat{\psi}^n_{i,j})\widehat{\phi}^{n+\frac 12}_{i,j} - f(x_j,\tilde{\phi}^{n+\frac 12}_{i,j}\tilde{\psi}^n_{i,j})\tilde{\phi}^{n+\frac 12}_{i,j}\right) \\
   &=& \delta^{n+\frac 12}_{i,j} \left(f(x_j,\widehat{\phi}^{n+\frac 12}_{i,j}\widehat{\psi}^n_{i,j})\widehat{\phi}^{n+\frac 12}_{i,j} - f(x_j,\widehat{\phi}^{n+\frac 12}_{i,j}\tilde{\psi}^n_{i,j})\widehat{\phi}^{n+\frac 12}_{i,j}\right) \\
  &+&  \left(\widehat{\phi}^{n + \frac 12}_{i,j} -  \tilde{\phi}^{n + \frac 12}_{i,j}\right)\left(f(x_j,\widehat{\phi}^{n+\frac 12}_{i,j}\tilde{\psi}^n_{i,j})\widehat{\phi}^{n+\frac 12}_{i,j} - f(x_j,\tilde{\phi}^{n+\frac 12}_{i,j}\tilde{\psi}^n_{i,j})\tilde{\phi}^{n+\frac 12}_{i,j}\right)\\
  &\le& \delta^{n+\frac 12}_{i,j} \left(f(x_j,\widehat{\phi}^{n+\frac 12}_{i,j}\widehat{\psi}^n_{i,j})\widehat{\phi}^{n+\frac 12}_{i,j} - f(x_j,\widehat{\phi}^{n+\frac 12}_{i,j}\tilde{\psi}^n_{i,j})\widehat{\phi}^{n+\frac 12}_{i,j}\right) \\
  &\le& K\left|\delta^{n+\frac 12}_{i,j}\right| \left|\widehat{\phi}^{n+\frac 12}_{i,j}\right| \left|\widehat{\phi}^{n+\frac 12}_{i,j}\widehat{\psi}^n_{i,j} - \widehat{\phi}^{n+\frac 12}_{i,j}\tilde{\psi}^n_{i,j}\right| \\
  &\le& K\left|\delta^{n+\frac 12}_{i,j}\right| \left|\widehat{\phi}^{n+\frac 12}_{i,j}\right|^2 \left|\delta^n_{i,j}\right| \\
\end{eqnarray*}
\begin{eqnarray*}
  &\le& \frac K2 \left\|\widehat{\phi}^{n+\frac 12}\right\|_{\infty}^2 \left[ \left(\delta^{n+\frac 12}_{i,j}\right)^2 + \left(\delta^n_{i,j}\right)^2\right] \\
  &\le& \frac K2 \left\|\widehat{\phi}^{\frac 12}\right\|_{\infty}^2 \left[ \left(\delta^{n+\frac 12}_{i,j}\right)^2 + \left(\delta^n_{i,j}\right)^2\right] \\
\end{eqnarray*}

Thus:

\begin{eqnarray*}
   & & \sum_{j=0}^J\left(\delta^{n+\frac 12}_{i,j}\right)^2 \left(1 -\Delta t - \frac{\Delta t}{\sigma^2} K \left\|\widehat{\phi}^{\frac 12}\right\|_{\infty}^2 \right) \\
   &\le&  \sum_{j=0}^J \left(\delta^{n+\frac 12}_{i+1,j}\right)^2 + \frac{\sigma^2 \Delta t}{ (\Delta x)^2} \sum_{j=0}^J \delta^{n+\frac 12}_{i,j} \left(\delta^{n+\frac 12}_{i,j+1} - 2\delta^{n+\frac 12}_{i,j} + \delta^{n+\frac 12}_{i,j-1}\right)\\
   &+& \frac{\Delta t}{\sigma^2} K \left\|\widehat{\phi}^{\frac 12}\right\|_{\infty}^2  \sum_{j=0}^J \left(\delta^{n}_{i,j}\right)^2 + \Delta t \sum_{j=0}^J \left(\tilde{\eta}^{n+\frac 12}_{i,j}\right)^2\\
\end{eqnarray*}
\begin{eqnarray*}
   &\le&  \sum_{j=0}^J \left(\delta^{n+\frac 12}_{i+1,j}\right)^2 - \frac{\sigma^2 \Delta t}{ (\Delta x)^2} \sum_{j=0}^{J-1} \left(\delta^{n+\frac 12}_{i,j+1} - \delta^{n+\frac 12}_{i,j}\right)^2\\
   &+& \frac{\Delta t}{\sigma^2} K \left\|\widehat{\phi}^{\frac 12}\right\|_{\infty}^2  (J+1) \left\|\delta^{n}\right\|^2 + \Delta t (J+1)\left\|\tilde{\eta}^{n+\frac 12}\right\|^2\\
   &\le&  \sum_{j=0}^J \left(\delta^{n+\frac 12}_{i+1,j}\right)^2 + \frac{\Delta t}{\sigma^2} K \left\|\widehat{\phi}^{\frac 12}\right\|_{\infty}^2  (J+1) \left\|\delta^{n}\right\|^2 + \Delta t (J+1)\left\|\tilde{\eta}^{n+\frac 12}\right\|^2\\
\end{eqnarray*}

We know that $\left\|\widehat{\phi}^{\frac 12}\right\|_{\infty} \le \left\|e^{\frac{u_T}{\sigma^2}}\right\|_{\infty}$. Hence, $\alpha(\Delta t) = 1 -\Delta t - \frac{\Delta t}{\sigma^2} K \left\|\widehat{\phi}^{\frac 12}\right\|_{\infty}^2$ lies in $(0,1)$ by hypothesis.\\

This gives:

$$ \sum_{j=0}^J \left(\delta^{n+\frac 12}_{i,j}\right)^2 \le \frac{1}{\alpha(\Delta t)}\sum_{j=0}^J \left(\delta^{n+\frac 12}_{i+1,j}\right)^2 + \frac{\Delta t}{\alpha(\Delta t)} (J+1)  \left[ \frac{1}{\sigma^2} K \left\|\widehat{\phi}^{\frac 12}\right\|_{\infty}^2  \left\|\delta^{n}\right\|^2 + \left\|\tilde{\eta}^{n+\frac 12}\right\|^2\right]$$

By immediate induction we have:

\begin{eqnarray*}
\sum_{j=0}^J \left(\delta^{n+\frac 12}_{i,j}\right)^2 &\le& \left(\frac{1}{\alpha(\Delta t)}\right)^{I-i}\sum_{j=0}^J \left(\underbrace{\delta^{n+\frac 12}_{I,j}}_{=0}\right)^2\\
&+& \Delta t \sum_{l=i}^{I-1} \left(\frac{1}{\alpha(\Delta t)}\right)^{l-i} (J+1)  \left[ \frac{1}{\sigma^2} K \left\|\widehat{\phi}^{\frac 12}\right\|_{\infty}^2  \left\|\delta^{n}\right\|^2 + \left\|\tilde{\eta}^{n+\frac 12}\right\|^2\right]\\
&\le& \left(\frac{1}{\alpha(\Delta t)}\right)^{I} T (J+1) \left[ \frac{1}{\sigma^2} K \left\|\widehat{\phi}^{\frac 12}\right\|_{\infty}^2  \left\|\delta^{n}\right\|^2 + \left\|\tilde{\eta}^{n+\frac 12}\right\|^2\right]\\
&\le& \left(\frac{1}{\alpha(\Delta t)}\right)^{\frac{T}{\Delta t}} T (J+1) \left[ \frac{1}{\sigma^2} K \left\|\widehat{\phi}^{\frac 12}\right\|_{\infty}^2  \left\|\delta^{n}\right\|^2 + \left\|\tilde{\eta}^{n+\frac 12}\right\|^2\right]\\
\end{eqnarray*}

Thus:

$$\left\|\delta^{n+\frac 12}\right\|^2 \le \left(\frac{1}{\alpha(\Delta t)}\right)^{\frac{T}{\Delta t}} T \left[ \frac{1}{\sigma^2} K \left\|\widehat{\phi}^{\frac 12}\right\|_{\infty}^2  \left\|\delta^{n}\right\|^2 + \left\|\tilde{\eta}^{n+\frac 12}\right\|^2\right]$$

Now, $\lim_{\Delta t \to 0}\left(\frac{1}{\alpha(\Delta t)}\right)^{\frac{T}{\Delta t}} = \exp\left(\left(1 + \frac{1}{\sigma^2} K \left\|\widehat{\phi}^{\frac 12}\right\|_{\infty}^2\right)T\right)$ and hence there exists a constant $C_1$ such that:

$$\left\|\delta^{n+\frac 12}\right\|^2 \le  C_1  \left\|\delta^{n}\right\|^2 + C_1 \left\|\tilde{\eta}^{n+\frac 12}\right\|^2$$

This is the first part of the Theorem. For the second part, we need to control the difference arising from the initial condition.\\

By definition we have for $i\in \mathcal{I}-\lbrace I \rbrace, j\in \mathcal{J}$:

$$\frac{\delta^{n+1}_{i+1,j} - \delta^{n+1}_{i,j}}{\Delta t} - \frac{\sigma^2}{2} \frac{\delta^{n+1}_{i+1,j+1} - 2\delta^{n+1}_{i+1,j} + \delta^{n+1}_{i+1,j-1}}{(\Delta x)^2}$$$$ = \frac{1}{\sigma^2} \left(f(x_j,\widehat{\phi}^{n+\frac 12}_{i+1,j}\widehat{\psi}^{n+1}_{i+1,j})\widehat{\psi}^{n+1}_{i+1,j} - f(x_j,\tilde{\phi}^{n+\frac 12}_{i+1,j}\tilde{\psi}^{n+1}_{i+1,j})\tilde{\psi}^{n+1}_{i+1,j}\right) - \tilde{\eta}^{n + 1}_{i+1,j}$$

Hence:

$$\left(\delta^{n+1}_{i+1,j}\right)^2 = \delta^{n+1}_{i+1,j} \delta^{n+1}_{i,j} + \frac{\sigma^2 \Delta t}{2 (\Delta x)^2} \delta^{n+1}_{i+1,j} \left(\delta^{n+1}_{i+1,j+1} - 2\delta^{n+1}_{i+1,j} + \delta^{n+1}_{i+1,j-1}\right)$$
$$+ \frac{\Delta t}{\sigma^2} \delta^{n+1}_{i+1,j} \left(f(x_j,\widehat{\phi}^{n+\frac 12}_{i+1,j}\widehat{\psi}^{n+1}_{i+1,j})\widehat{\psi}^{n+1}_{i+1,j} - f(x_j,\tilde{\phi}^{n+\frac 12}_{i+1,j}\tilde{\psi}^{n+1}_{i+1,j})\tilde{\psi}^{n+1}_{i+1,j}\right) - \delta^{n+1}_{i+1,j} \tilde{\eta}^{n+1}_{i+1,j} \Delta t$$
$$\Rightarrow \frac 12 \left(\delta^{n+1}_{i+1,j}\right)^2 \le \frac 12 \left(\delta^{n+1}_{i+1,j}\right)^2 + \frac{\sigma^2 \Delta t}{2 (\Delta x)^2} \delta^{n+1}_{i+1,j} \left(\delta^{n+1}_{i+1,j+1} - 2\delta^{n+1}_{i+1,j} + \delta^{n+1}_{i+1,j-1}\right)$$
$$+ \frac{\Delta t}{\sigma^2} \delta^{n+1}_{i+1,j} \left(f(x_j,\widehat{\phi}^{n+\frac 12}_{i+1,j}\widehat{\psi}^{n+1}_{i+1,j})\widehat{\psi}^{n+1}_{i+1,j} - f(x_j,\tilde{\phi}^{n+\frac 12}_{i+1,j}\tilde{\psi}^{n+1}_{i+1,j})\tilde{\psi}^{n+1}_{i+1,j}\right) + \frac {\Delta t}{2} \left(\delta^{n+1}_{i+1,j}\right)^2 + \frac {\Delta t}{2} \left(\tilde{\eta}^{n+1}_{i+1,j}\right)^2$$

Now,

\begin{eqnarray*}
   & & \delta^{n+1}_{i+1,j} \left(f(x_j,\widehat{\phi}^{n+\frac 12}_{i+1,j}\widehat{\psi}^{n+1}_{i+1,j})\widehat{\psi}^{n+1}_{i+1,j} - f(x_j,\tilde{\phi}^{n+\frac 12}_{i+1,j}\tilde{\psi}^{n+1}_{i+1,j})\tilde{\psi}^{n+1}_{i+1,j}\right) \\
   &=& \delta^{n+1}_{i+1,j} \left(f(x_j,\widehat{\phi}^{n+\frac 12}_{i+1,j}\widehat{\psi}^{n+1}_{i+1,j})\widehat{\psi}^{n+1}_{i+1,j} - f(x_j,\tilde{\phi}^{n+\frac 12}_{i+1,j}\widehat{\psi}^{n+1}_{i+1,j})\widehat{\psi}^{n+1}_{i+1,j}\right) \\
  &+&  \left(\widehat{\psi}^{n + 1}_{i,j} -  \tilde{\psi}^{n + 1}_{i,j}\right)\left(f(x_j,\tilde{\phi}^{n+\frac 12}_{i+1,j}\widehat{\psi}^{n+1}_{i+1,j})\widehat{\psi}^{n+1}_{i+1,j} - f(x_j,\tilde{\phi}^{n+\frac 12}_{i+1,j}\tilde{\psi}^{n+1}_{i+1,j})\tilde{\psi}^{n+1}_{i+1,j}\right)\\
  &\le& \delta^{n+1}_{i+1,j} \left(f(x_j,\widehat{\phi}^{n+\frac 12}_{i+1,j}\widehat{\psi}^{n+1}_{i+1,j})\widehat{\psi}^{n+1}_{i+1,j} - f(x_j,\tilde{\phi}^{n+\frac 12}_{i+1,j}\widehat{\psi}^{n+1}_{i+1,j})\widehat{\psi}^{n+1}_{i+1,j}\right) \\
  &\le& K\left|\delta^{n+1}_{i+1,j}\right| \left|\widehat{\psi}^{n+1}_{i+1,j}\right| \left|\widehat{\phi}^{n+\frac 12}_{i+1,j}\widehat{\psi}^{n+1}_{i+1,j} - \tilde{\phi}^{n+\frac 12}_{i+1,j}\widehat{\psi}^{n+1}_{i+1,j}\right| \\
  &\le& K\left|\delta^{n+1}_{i+1,j}\right| \left|\widehat{\psi}^{n+1}_{i+1,j}\right|^2 \left|\delta^{n+\frac 12}_{i+1,j}\right| \\
  &\le& \frac K2 \left\|\widehat{\psi}^{n+1}\right\|_{\infty}^2 \left[ \left(\delta^{n+\frac 12}_{i+1,j}\right)^2 + \left(\delta^{n+1}_{i+1,j}\right)^2\right] \\
  &\le& \frac K2 \left\|\widehat{\psi}\right\|_{\infty}^2 \left[ \left(\delta^{n+\frac 12}_{i+1,j}\right)^2 + \left(\delta^{n+1}_{i+1,j}\right)^2\right] \\
\end{eqnarray*}

Thus:

\begin{eqnarray*}
   & & \sum_{j=0}^J\left(\delta^{n+1}_{i+1,j}\right)^2 \left(1 -\Delta t - \frac{\Delta t}{\sigma^2} K \left\|\widehat{\psi}\right\|_{\infty}^2 \right) \\
   &\le&  \sum_{j=0}^J \left(\delta^{n+1}_{i,j}\right)^2 + \frac{\sigma^2 \Delta t}{ (\Delta x)^2} \sum_{j=0}^J \delta^{n+1}_{i+1,j} \left(\delta^{n+1}_{i+1,j+1} - 2\delta^{n+1}_{i+1,j} + \delta^{n+1}_{i+1,j-1}\right)\\
   &+& \frac{\Delta t}{\sigma^2} K \left\|\widehat{\psi}\right\|_{\infty}^2  \sum_{j=0}^J \left(\delta^{n+\frac 12}_{i+1,j}\right)^2 + \Delta t \sum_{j=0}^J \left(\tilde{\eta}^{n+1}_{i+1,j}\right)^2\\
   &\le&  \sum_{j=0}^J \left(\delta^{n+1}_{i,j}\right)^2 - \frac{\sigma^2 \Delta t}{ (\Delta x)^2} \sum_{j=0}^{J-1} \left(\delta^{n+1}_{i+1,j+1} - \delta^{n+1}_{i+1,j}\right)^2\\
   &+& \frac{\Delta t}{\sigma^2} K \left\|\widehat{\psi}\right\|_{\infty}^2  (J+1) \left\|\delta^{n+\frac 12}\right\|^2 + \Delta t (J+1)\left\|\tilde{\eta}^{n+1}\right\|^2\\
   &\le&  \sum_{j=0}^J \left(\delta^{n+1}_{i,j}\right)^2 + \frac{\Delta t}{\sigma^2} K \left\|\widehat{\psi}\right\|_{\infty}^2  (J+1) \left\|\delta^{n+\frac{1}{2}}\right\|^2 + \Delta t (J+1)\left\|\tilde{\eta}^{n+1}\right\|^2\\
\end{eqnarray*}

Let us introduce $\beta(\Delta t) = 1 -\Delta t - \frac{\Delta t}{\sigma^2} K \left\|\widehat{\psi}\right\|_{\infty}^2$. By hypothesis, $\beta$ lies in $(0,1)$.\\

This gives:

$$ \beta(\Delta t) \sum_{j=0}^J \left(\delta^{n+1}_{i+1,j}\right)^2 \le \sum_{j=0}^J \left(\delta^{n+1}_{i,j}\right)^2 + \Delta t (J+1)  \left[ \frac{1}{\sigma^2} K \left\|\widehat{\psi}\right\|_{\infty}^2  \left\|\delta^{n+\frac 12}\right\|^2 + \left\|\tilde{\eta}^{n+1}\right\|^2\right]$$

By immediate induction we have:

\begin{eqnarray*}
\left(\beta(\Delta t)\right)^{i+1} \sum_{j=0}^J \left(\delta^{n+1}_{i+1,j}\right)^2 &\le& \sum_{j=0}^J \left(\delta^{n+1}_{0,j}\right)^2 +  \Delta t \sum_{l=0}^{i} \left(\beta(\Delta t)\right)^{l} (J+1)  \left[ \frac{1}{\sigma^2} K \left\|\widehat{\psi}\right\|_{\infty}^2  \left\|\delta^{n+\frac 12}\right\|^2 + \left\|\tilde{\eta}^{n+1}\right\|^2\right]\\
\end{eqnarray*}

Hence:

$$\sum_{j=0}^J \left(\delta^{n+1}_{i+1,j}\right)^2 \le \left(\frac{1}{\beta(\Delta t)}\right)^{I} \sum_{j=0}^J \left(\delta^{n+1}_{0,j}\right)^2$$$$+  \left(\frac{1}{\beta(\Delta t)}\right)^{I} T (J+1) \left[ \frac{1}{\sigma^2} K \left\|\widehat{\psi}\right\|_{\infty}^2  \left\|\delta^{n+\frac 12}\right\|^2 + \left\|\tilde{\eta}^{n+1}\right\|^2\right]$$

Now, $$\delta^{n+1}_{0,j} = \frac{m_0(x_j)}{\widehat{\phi}^{n+\frac 12}_{0,j}} - \frac{m_0(x_j)}{\tilde{\phi}^{n+\frac 12}_{0,j}}$$ Consequently:
$$|\delta^{n+1}_{0,j}| \le \|m_0\|_{\infty} \left|\frac{1}{\widehat{\phi}^{n+\frac 12}_{0,j}} - \frac{1}{\tilde{\phi}^{n+\frac 12}_{0,j}}\right| \le \frac{\|m_0\|_{\infty}}{\epsilon^2} \left|\delta^{n+\frac 12}_{0,j}\right|$$

We end up with:

$$\sum_{j=0}^J \left(\delta^{n+1}_{i+1,j}\right)^2 \le \left(\frac{1}{\beta(\Delta t)}\right)^{I} \left(\frac{\|m_0\|_{\infty}}{\epsilon^2}\right)^2 (J+1)\left\|\delta^{n+\frac 12}\right\|^2 $$$$+  \left(\frac{1}{\beta(\Delta t)}\right)^{I} T (J+1) \left[ \frac{1}{\sigma^2} K \left\|\widehat{\psi}\right\|_{\infty}^2  \left\|\delta^{n+\frac 12}\right\|^2 + \left\|\tilde{\eta}^{n+1}\right\|^2\right]$$

$$\Rightarrow \left\|\delta^{n+1}\right\|^2 \le \left(\frac{1}{\beta(\Delta t)}\right)^{I} \left(\frac{\|m_0\|_{\infty}}{\epsilon^2}\right)^2 (J+1)\left\|\delta^{n+\frac 12}\right\|^2 $$$$+  \left(\frac{1}{\beta(\Delta t)}\right)^{I} T (J+1) \left[ \frac{1}{\sigma^2} K \left\|\widehat{\psi}\right\|_{\infty}^2  \left\|\delta^{n+\frac 12}\right\|^2 + \left\|\tilde{\eta}^{n+1}\right\|^2\right]$$

Now, as above, $\lim_{\Delta t \to 0}\left(\frac{1}{\beta(\Delta t)}\right)^{\frac{T}{\Delta t}} = \exp\left(\left(1 + \frac{1}{\sigma^2} K \left\|\widehat{\psi}\right\|_{\infty}^2\right)T\right)$ and hence there exists a constant $C_2$ such that:

$$\left\|\delta^{n+1}\right\|^2 \le  C_2  \left\|\delta^{n+\frac 12}\right\|^2 + C_2 \left\|\tilde{\eta}^{n+1}\right\|^2$$

Eventually, $C=\max(C_1,C_2)$ provides the bounds given in the theorem\qed\\

These bounds allow to prove, under regularity assumptions, that the scheme is convergent. Let us start with the convergence of the schemes defining respectively $\widehat{\phi}^{n+\frac 12}$ and $\widehat{\psi}^{n}$.

\begin{Theorem}
Let us suppose, in addition to the hypotheses of Theorem 2 that $u_T$, $m_0$ and $f$ are so that $\forall n \in \mathbb{N}, \phi^{n+\frac{1}{2}}, \psi^n \in C^{1,2}([0,T]\times[0,1])$.\\
Then $\forall n \in \mathbb{N}$:
$$\lim_{\Delta t, \Delta x \to 0} \|\widehat{\phi}^{n+\frac 12} - \tilde{\phi}^{n+\frac 12}\| =0$$
$$ \lim_{\Delta t, \Delta x \to 0} \|\widehat{\psi}^{n+1} - \tilde{\psi}^{n+1}\| =0$$
\end{Theorem}

\noindent\textbf{Proof:}\\

Let us fix $n\ge 0$. We have by immediate induction from Theorem 2 that, as soon as $\Delta t$ is small enough, that there exists a constant $C$ independent of $\Delta t$ and $\Delta x$ so that:

$$\|\widehat{\phi}^{n+\frac 12} - \tilde{\phi}^{n+\frac 12}\|^2 \le \sum_{k=0}^{n} C^{2k+1} \|\tilde{\eta}^{n+\frac 12 -k}\|^2 + C^{2k+2} \|\tilde{\eta}^{n-k}\|^2$$
and
$$\|\widehat{\psi}^{n+1} - \tilde{\psi}^{n+1}\|^2 \le \sum_{k=0}^{n} C^{2k+1} \|\tilde{\eta}^{n+1-k}\|^2 + C^{2k+2} \|\tilde{\eta}^{n+\frac 12-k}\|^2$$

Hence, because our assumptions imply that the consistency errors tend to $0$ as $\Delta t$ and $\Delta x$ tend to $0$, we have that:

$$\lim_{\Delta t, \Delta x \to 0} \|\widehat{\phi}^{n+\frac 12} - \tilde{\phi}^{n+\frac 12}\| = 0$$
and
$$\lim_{\Delta t, \Delta x \to 0} \|\widehat{\psi}^{n+1} - \tilde{\phi}^{n+1}\| = 0$$
\qed\\

For the global convergence of the numerical scheme we have the following theorem:

\begin{Theorem}
Let us suppose, in addition to the hypotheses made in Theorem 2 that $u_T$, $m_0$ and $f$ are so that $\forall n \in \mathbb{N}, \phi^{n+\frac{1}{2}}, \psi^n \in C^{1,2}([0,T]\times[0,1])$ and $\phi, \psi \in C^{1,2}([0,T]\times[0,1])$\\
Then:
$$\lim_{n\to \infty} \limsup_{\Delta t, \Delta x \to 0} \|\widehat{\phi}^{n+\frac 12} - \tilde{\phi}\| =0$$
$$\lim_{n\to \infty} \limsup_{\Delta t, \Delta x \to 0} \|\widehat{\psi}^{n+1} - \tilde{\psi}\| =0$$
\end{Theorem}

\noindent\textbf{Proof:}\\

We have seen that $(\phi^{n+\frac 12})_n$ converges monotonically towards $\phi$. Hence, we know from Dini's theorem that the convergence is in fact uniform. Consequently:
$$\forall \varepsilon >0, \exists n_0, \forall n\ge n_0, \forall \Delta t, \Delta x > 0, \|\tilde{\phi}^{n+\frac 12} - \tilde{\phi}\| \le \|\tilde{\phi}^{n+\frac 12} - \tilde{\phi}\|_{\infty} \le \varepsilon$$
Hence, $\forall n\ge n_0$, we know from the preceding result that:
$$\limsup_{\Delta t, \Delta x \to 0} \|\widehat{\phi}^{n+\frac 12} - \tilde{\phi}\| \le  \limsup_{\Delta t, \Delta x \to 0} \|\widehat{\phi}^{n+\frac 12} - \tilde{\phi}^{n+\frac 12}\|+\|\tilde{\phi}^{n+\frac 12} - \tilde{\phi}\| \le \varepsilon$$

Thus, we get that:

$$\lim_{n\to \infty} \limsup_{\Delta t, \Delta x \to 0} \|\widehat{\phi}^{n+\frac 12} - \tilde{\phi}\| =0$$

The same proof holds regarding $\tilde{\psi}$.\qed\\

We have proved that the scheme was convergent and we provided a bound on $\Delta t$ that guarantees stability of the scheme. Although this bound is certainly not the best one, a remark must be done regarding how it depends on parameters and especially on $\sigma$.\\
It's indeed important to notice that the bound in Theorem 2 is $$1+\frac{K}{\sigma^2}\max\left(\left\|e^{\frac{u_T}{\sigma^2}} \right\|^2_\infty,\frac{\|m_0\|^2_\infty}{\epsilon^2}\right)=1+\frac{K}{\sigma^2}\max\left(\left\|e^{\frac{u_T}{\sigma^2}} \right\|^2_\infty,\|m_0\|^2_\infty \exp\left(\frac{2}{\sigma^2}\left(\|u_T\|_{\infty} + \|f\|_{\infty} T\right)\right)\right)$$ and this bound tends very rapidly to $+\infty$ as $\sigma$ tends to $0$. As a consequence, the numerical scheme may not be useful for small values of $\sigma$.\\

We will empirically tackle the question of the influence of $\sigma$ in the next section but before turning to the numerical experiments, let's highlight another important property of this numerical scheme.\\
In most cases in which a probability distribution function is involved, the numerical scheme is built so that mass is preserved along the trajectory. Here, as we have seen above, since $m$ is only reconstructed at the end from $\psi$ and $\phi$, there is nothing like mass conservation. However, a difference arises with Proposition 6 since there is theoretically no systematic loss of mass for a fixed $n$ as $t$ goes from $0$ to $T$. The evolution of total mass is given by the following proposition:\\

\begin{Proposition}
Let us introduce $\widehat{m}^{n+1} = \widehat{\phi}^{n+\frac 12} \widehat{\psi}^{n+1}$.\\
The difference in total mass between two successive times at step $n+1$ is:
$$\frac{1}{J+1} \sum_{j=0}^J \widehat{m}_{i+1,j}^{n+1} - \frac{1}{J+1} \sum_{j=0}^J \widehat{m}_{i,j}^{n+1} = \frac{\Delta t}{\sigma^2}\frac{1}{J+1} \sum_{j=0}^J \widehat{\psi}_{i+1,j}^{n+1} \widehat{\phi}_{i,j}^{n+\frac 12}\left(f\left(\widehat{\psi}_{i+1,j}^{n+1} \widehat{\phi}_{i+1,j}^{n+\frac 12}\right)-f\left(\widehat{\psi}_{i,j}^{n} \widehat{\phi}_{i,j}^{n+\frac 12}\right)\right)$$
Then, under the hypotheses of Theorem 4, we have that:
$$\forall i\in \mathcal{I}, \lim_{n\to\infty} \limsup_{\Delta t,\Delta x \to 0} \frac{1}{J+1} \sum_{j=0}^J \widehat{m}_{i,j}^{n+1} = \int_{\Omega} m_0(x) dx$$
\end{Proposition}

\noindent \textbf{Proof:}\\

\begin{eqnarray*}
   & & \frac{1}{J+1} \sum_{j=0}^J \widehat{m}_{i+1,j}^{n+1} - \frac{1}{J+1} \sum_{j=0}^J \widehat{m}_{i,j}^{n+1} \\
    &=&  \frac{1}{J+1} \sum_{j=0}^J \widehat{\psi}_{i+1,j}^{n+1}\widehat{\phi}_{i+1,j}^{n+\frac 12} - \frac{1}{J+1} \sum_{j=0}^J \widehat{\psi}_{i,j}^{n+1}\widehat{\phi}_{i,j}^{n+\frac 12} \\
\end{eqnarray*}
\begin{eqnarray*}
   &=&  \frac{1}{J+1} \sum_{j=0}^J \widehat{\psi}_{i+1,j}^{n+1}\left(\widehat{\phi}_{i+1,j}^{n+\frac 12} - \widehat{\phi}_{i,j}^{n+\frac 12} \right) + \frac{1}{J+1} \sum_{j=0}^J \left(\widehat{\psi}_{i+1,j}^{n+1} - \widehat{\psi}_{i,j}^{n+1}\right)\widehat{\phi}_{i,j}^{n+\frac 12} \\
   &=&  -\frac{\Delta t}{J+1} \left[ \frac{\sigma^2}{2(\Delta x)^2} \sum_{j=0}^J \widehat{\psi}_{i+1,j}^{n+1}\left(\widehat{\phi}_{i,j+1}^{n+\frac 12} - 2 \widehat{\phi}_{i,j}^{n+\frac 12} + \widehat{\phi}_{i,j-1}^{n+\frac 12} \right) + \frac{1}{\sigma^2} \sum_{j=0}^J \widehat{\psi}_{i+1,j}^{n+1}
   \widehat{\phi}_{i,j}^{n+\frac 12} f\left(\widehat{\phi}_{i,j}^{n+\frac 12}\widehat{\psi}_{i,j}^{n}\right)\right]\\
   &+&  \frac{\Delta t}{J+1} \left[\frac{\sigma^2}{2(\Delta x)^2} \sum_{j=0}^J \widehat{\phi}_{i,j}^{n+\frac 12}\left(\widehat{\psi}_{i+1,j+1}^{n+1} - 2 \widehat{\psi}_{i+1,j}^{n+1} + \widehat{\psi}_{i+1,j-1}^{n+1} \right) + \frac{1}{\sigma^2} \sum_{j=0}^J \widehat{\psi}_{i+1,j}^{n+1}
   \widehat{\phi}_{i,j}^{n+\frac 12} f\left(\widehat{\phi}_{i+1,j}^{n+\frac 12}\widehat{\psi}_{i+1,j}^{n+1}\right)\right]\\
   &=&  -\frac{\Delta t}{J+1} \left[ -\frac{\sigma^2}{2(\Delta x)^2} \sum_{j=0}^{J-1} \left(\widehat{\psi}_{i+1,j+1}^{n+1} - \widehat{\psi}_{i+1,j}^{n+1} \right)\left(\widehat{\phi}_{i,j+1}^{n+\frac 12} - \widehat{\phi}_{i,j}^{n+\frac 12}\right) + \frac{1}{\sigma^2} \sum_{j=0}^J \widehat{\psi}_{i+1,j}^{n+1}\widehat{\phi}_{i,j}^{n+\frac 12} f\left(\widehat{\phi}_{i,j}^{n+\frac 12}\widehat{\psi}_{i,j}^{n}\right)\right]\\
   &+&  \frac{\Delta t}{J+1}\left[ -\frac{\sigma^2}{2(\Delta x)^2} \sum_{j=0}^{J-1} \left(\widehat{\psi}_{i+1,j+1}^{n+1} - \widehat{\psi}_{i+1,j}^{n+1} \right)\left(\widehat{\phi}_{i,j+1}^{n+\frac 12} - \widehat{\phi}_{i,j}^{n+\frac 12}\right) + \frac{1}{\sigma^2}\sum_{j=0}^J \widehat{\psi}_{i+1,j}^{n+1} \widehat{\phi}_{i,j}^{n+\frac 12} f\left(\widehat{\phi}_{i+1,j}^{n+\frac 12}\widehat{\psi}_{i+1,j}^{n+1}\right)\right]\\
   &=& \frac{\Delta t}{\sigma^2}\frac{1}{J+1} \sum_{j=0}^J \widehat{\psi}_{i+1,j}^{n+1} \widehat{\phi}_{i,j}^{n+\frac 12}\left(f\left(\widehat{\psi}_{i+1,j}^{n+1} \widehat{\phi}_{i+1,j}^{n+\frac 12}\right)-f\left(\widehat{\psi}_{i,j}^{n} \widehat{\phi}_{i,j}^{n+\frac 12}\right)\right)
\end{eqnarray*}

Hence, if the assumptions are the same as in Theorem 4, we have that:

\begin{eqnarray*}
   & & \left|\frac{1}{J+1} \sum_{j=0}^J \widehat{m}_{i+1,j}^{n+1} - \frac{1}{J+1} \sum_{j=0}^J \widehat{m}_{i,j}^{n+1}\right| \\
    &\le &  \frac{\Delta t}{\sigma^2}\frac{1}{J+1} K \left\|\widehat{\psi}^{n+1}\right\|_{\infty} \left\|\widehat{\phi}^{n+\frac 12}\right\|_{\infty} \sum_{j=0}^J \left|\widehat{\psi}_{i+1,j}^{n+1} \widehat{\phi}_{i+1,j}^{n+\frac 12}-\widehat{\psi}_{i,j}^{n} \widehat{\phi}_{i,j}^{n+\frac 12}\right|\\
    &\le &  \frac{\Delta t}{\sigma^2}\frac{1}{J+1} K \left\|e^{\frac{u_T}{\sigma^2}} \right\|_\infty \frac{\|m_0\|_\infty}{\epsilon}\sum_{j=0}^J \left|\widehat{\psi}_{i+1,j}^{n+1}\right| \left|\widehat{\phi}_{i+1,j}^{n+\frac 12} -\widehat{\phi}_{i,j}^{n+\frac 12}\right| +\left|\widehat{\phi}_{i,j}^{n+\frac 12}\right|\left|\widehat{\psi}_{i+1,j}^{n+1}-\widehat{\psi}_{i,j}^{n}\right| \\
    &\le &  \frac{\Delta t}{\sigma^2}\frac{1}{J+1} K \left\|e^{\frac{u_T}{\sigma^2}} \right\|_\infty \frac{\|m_0\|_\infty}{\epsilon} \max\left(\left\|e^{\frac{u_T}{\sigma^2}} \right\|_\infty,\frac{\|m_0\|_\infty}{\epsilon}\right)\sum_{j=0}^J \left|\widehat{\phi}_{i+1,j}^{n+\frac 12} -\widehat{\phi}_{i,j}^{n+\frac 12}\right| + \left|\widehat{\psi}_{i+1,j}^{n+1}-\widehat{\psi}_{i,j}^{n}\right| \\
\end{eqnarray*}

Now,

\begin{eqnarray*}
 & & \frac{1}{J+1}\sum_{j=0}^J \left|\widehat{\phi}_{i+1,j}^{n+\frac 12} -\widehat{\phi}_{i,j}^{n+\frac 12}\right|\\
 & \le & \frac{1}{J+1}\sum_{j=0}^J \left|\widehat{\phi}_{i+1,j}^{n+\frac 12} -\tilde{\phi}_{i+1,j}^{n+\frac 12}\right| + \left|\tilde{\phi}_{i+1,j}^{n+\frac 12} -\tilde{\phi}_{i,j}^{n+\frac 12}\right| + \left|\tilde{\phi}_{i,j}^{n+\frac 12} - \widehat{\phi}_{i,j}^{n+\frac 12}\right|\\
 & \le & 2 \left\|\widehat{\phi}^{n+\frac 12} -\tilde{\phi}^{n+\frac 12}\right\| + \frac{1}{J+1}\sum_{j=0}^J \left|\tilde{\phi}_{i+1,j}^{n+\frac 12} -\tilde{\phi}_{i,j}^{n+\frac 12}\right| \\
\end{eqnarray*}

and

\begin{eqnarray*}
 & & \frac{1}{J+1}\sum_{j=0}^J \left|\widehat{\psi}_{i+1,j}^{n+1} -\widehat{\psi}_{i,j}^{n}\right|\\
 & \le & \left\|\widehat{\psi}^{n+1} -\tilde{\psi}^{n+1}\right\| + \frac{1}{J+1}\sum_{j=0}^J \left|\tilde{\psi}_{i+1,j}^{n+1} -\tilde{\psi}_{i,j}^{n}\right| + \left\|\widehat{\psi}^{n} -\tilde{\psi}^{n}\right\| \\
\end{eqnarray*}

As a consequence, $\forall i\in \mathcal{I}$ :

\begin{eqnarray*}
   & & \left|\frac{1}{J+1} \sum_{j=0}^J \widehat{m}_{i,j}^{n+1} - \frac{1}{J+1} \sum_{j=0}^J \widehat{m}^{n+1}_{0,j}\right| \\
   & \le &  \frac{T}{\sigma^2}K \left\|e^{\frac{u_T}{\sigma^2}} \right\|_\infty \frac{\|m_0\|_\infty}{\epsilon} \max\left(\left\|e^{\frac{u_T}{\sigma^2}} \right\|_\infty,\frac{\|m_0\|_\infty}{\epsilon}\right)\\
   & \times &  \left[2 \left\|\widehat{\phi}^{n+\frac 12} -\tilde{\phi}^{n+\frac 12}\right\| + \left\|\widehat{\psi}^{n+1} -\tilde{\psi}^{n+1}\right\| + \left\|\widehat{\psi}^{n} -\tilde{\psi}^{n}\right\| + \frac{1}{J+1}\sup_{l\in \mathcal{I}} \sum_{j=0}^J  \left|\tilde{\phi}_{l+1,j}^{n+\frac 12} -\tilde{\phi}_{l,j}^{n+\frac 12}\right| + \left|\tilde{\psi}_{l+1,j}^{n+1} -\tilde{\psi}_{l,j}^{n}\right|\right]\\
\end{eqnarray*}

Using Theorem 4, we obtain that
$$\lim_{n\to\infty} \limsup_{\Delta t,\Delta x \to 0} \frac{1}{J+1} \sum_{j=0}^J \widehat{m}_{i,j}^{n+1} - \frac{1}{J+1} \sum_{j=0}^J \widehat{m}^{n+1}_{0,j} = 0$$

But, $\frac{1}{J+1} \sum_{j=0}^J \widehat{m}^{n+1}_{0,j} = \frac{1}{J+1} \sum_{j=0}^J m_0(x_j)$ and this expression tends to $\int_{\Omega} m_0(x) dx$ when $\Delta x$ tends to $0$ since $m_0$ is continuous from the hypotheses of Theorem 4.\qed\\
\section{Numerical examples}

Let us consider, as an example, the following problem on $\Omega = (0,1)$. Agents want to live near the center of the domain but prefer to live far from the others. This is modeled by the following\footnote{The bounds on $\xi$ guarantee that $f$ is bounded. In practice, these bounds are not binding.} function $f$:
$$f(x,\xi) = -16(x-1/2)^2-0.1\max(0,\min(5,\xi))$$

At the end of the period, when $t=T=0.5$, we suppose that agents do not have any incentive, namely $u_T(x)=0$.\\
We suppose that volatility is $\sigma=1$.\\

We assume that the population is initially distributed as: $$m_0(x) = \frac{\mu(x)}{\int_0^1 \mu(x')dx'}\quad  \mathrm{where} \quad \mu(x) = 1 + 0.2\cos\left(\pi\left(2x-\frac{3}{2}\right)\right)^2$$

In this case, we solved the numerical scheme with a Newton method for $\Delta t = 0.01$ and $\Delta x = 0.02$ and we obtained the following results for\footnote{$\alpha=\nabla u$ is the optimal control in this context.} $(\phi,\psi,u,m,\nabla u)$. The iterations in $n$ were stopped once the difference $\left\|\widehat{\phi}^{n+\frac 12}\widehat{\psi}^{n+1} - \widehat{\phi}^{n-\frac 12}\widehat{\psi}^{n} \right\|_{\infty}$ was below $10^{-7}$ (here 5 steps).\\

We see that quite rapidly the population as a whole gets close to what would be the stationary equilibrium. Then, just before the time horizon, the agents have no more incentive to be around $0$ and do not pay the cost associated to the stationary equilibrium. Hence the distribution spreads just before the end.\\

\begin{figure}[!h]
  \centering
  \includegraphics[width=420pt]{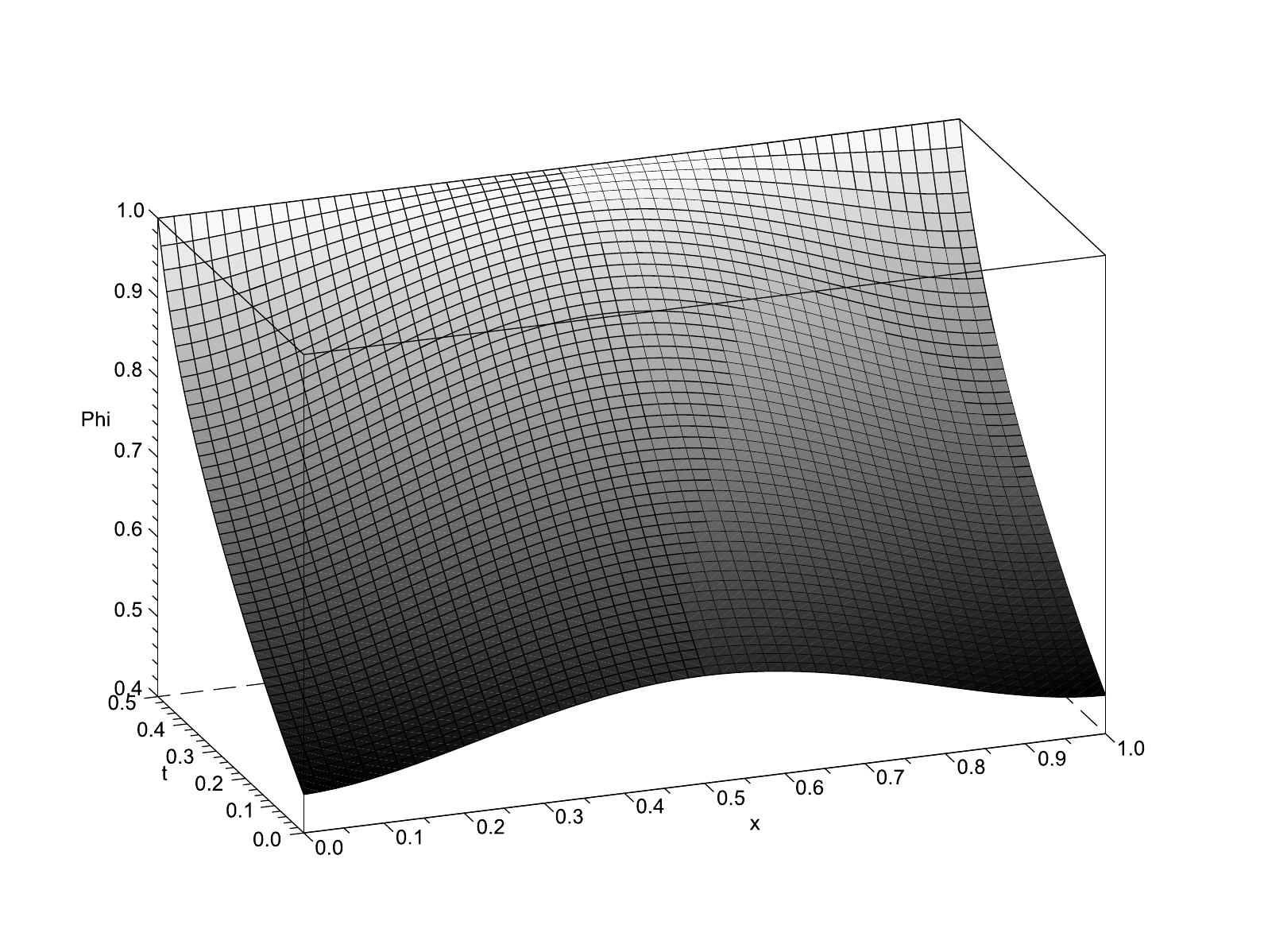}\\
  \caption{Solution for $\phi$}
\end{figure}

\begin{figure}[!h]
  \centering
  \includegraphics[width=420pt]{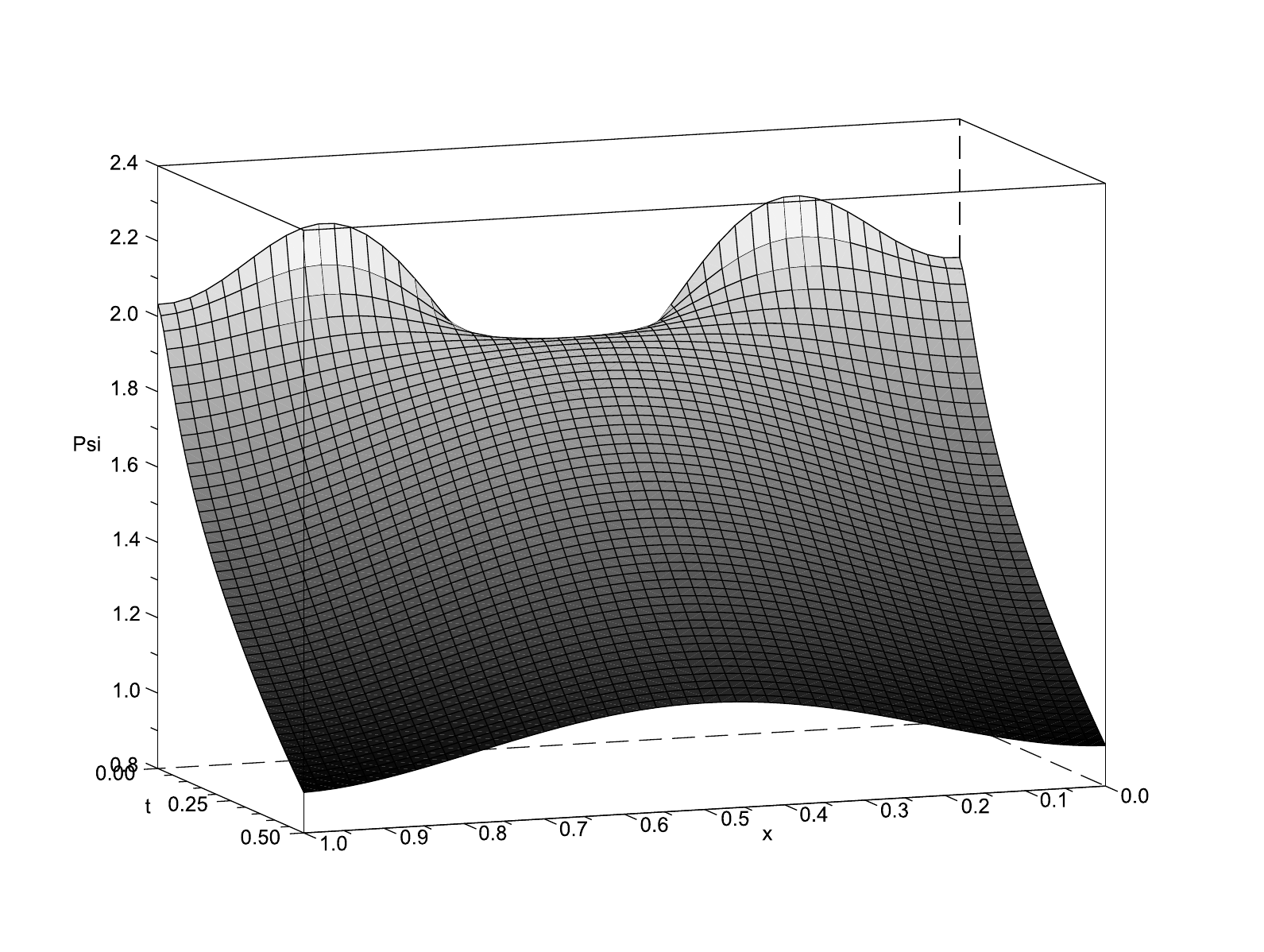}\\
  \caption{Solution for $\psi$}
\end{figure}
\newpage
\begin{figure}[!h]
  \centering
  \includegraphics[width=420pt]{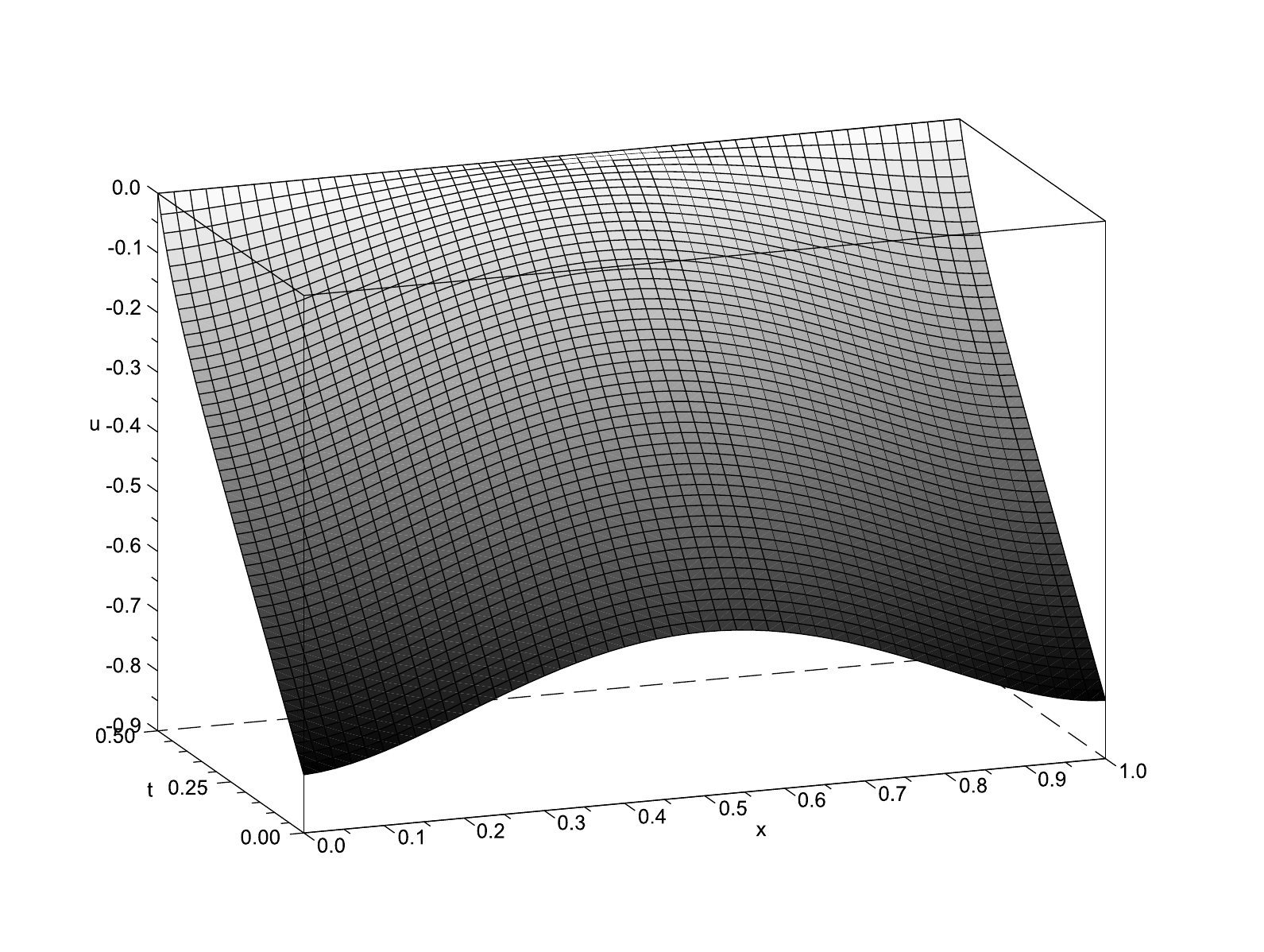}\\
  \caption{Solution for $u$}
\end{figure}

\begin{figure}[!h]
\centering
  \includegraphics[width=420pt]{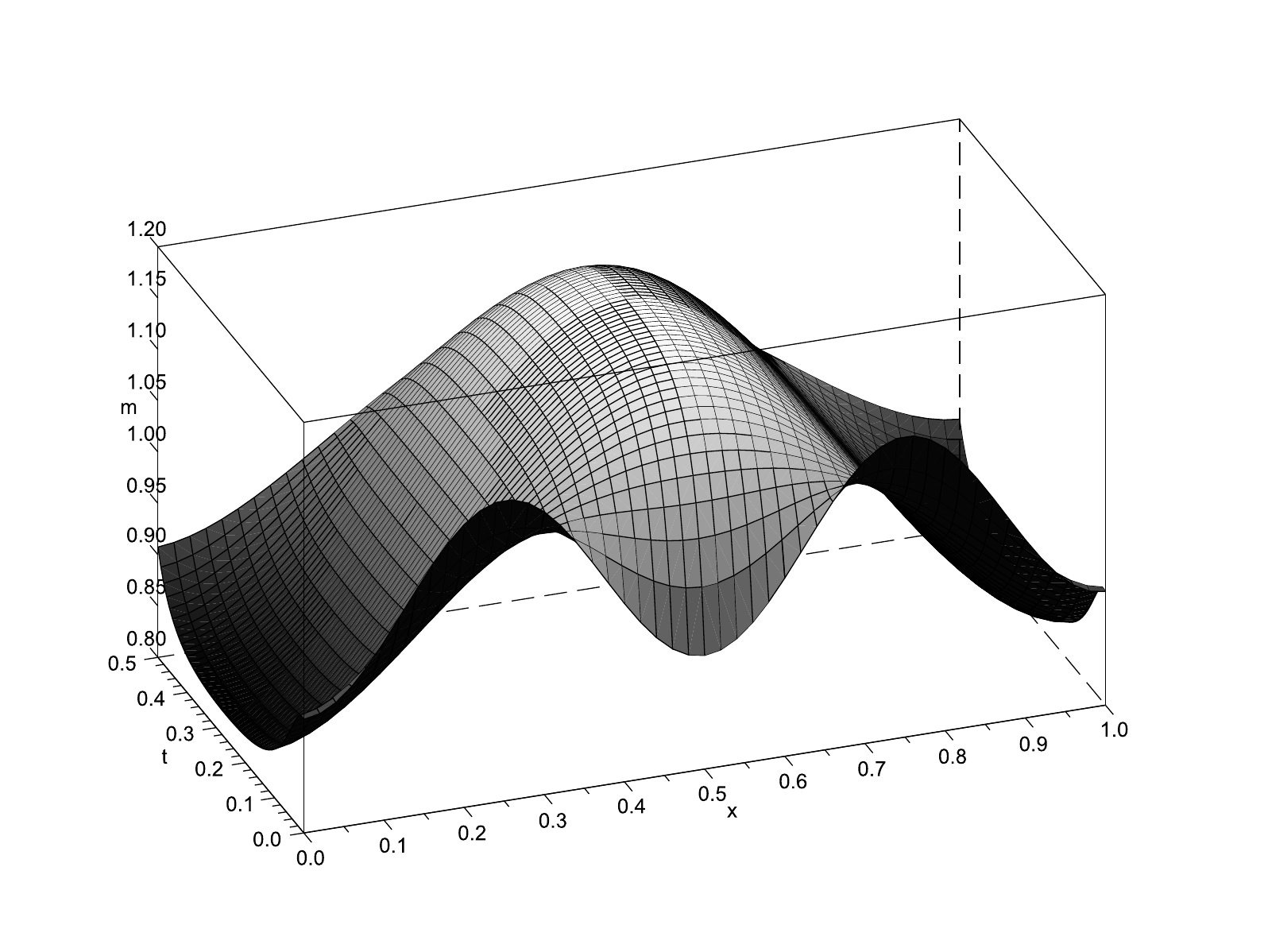}\\
  \caption{Solution for $m$}
\end{figure}

\newpage

\begin{figure}[!h]
\centering
  \includegraphics[width=420pt]{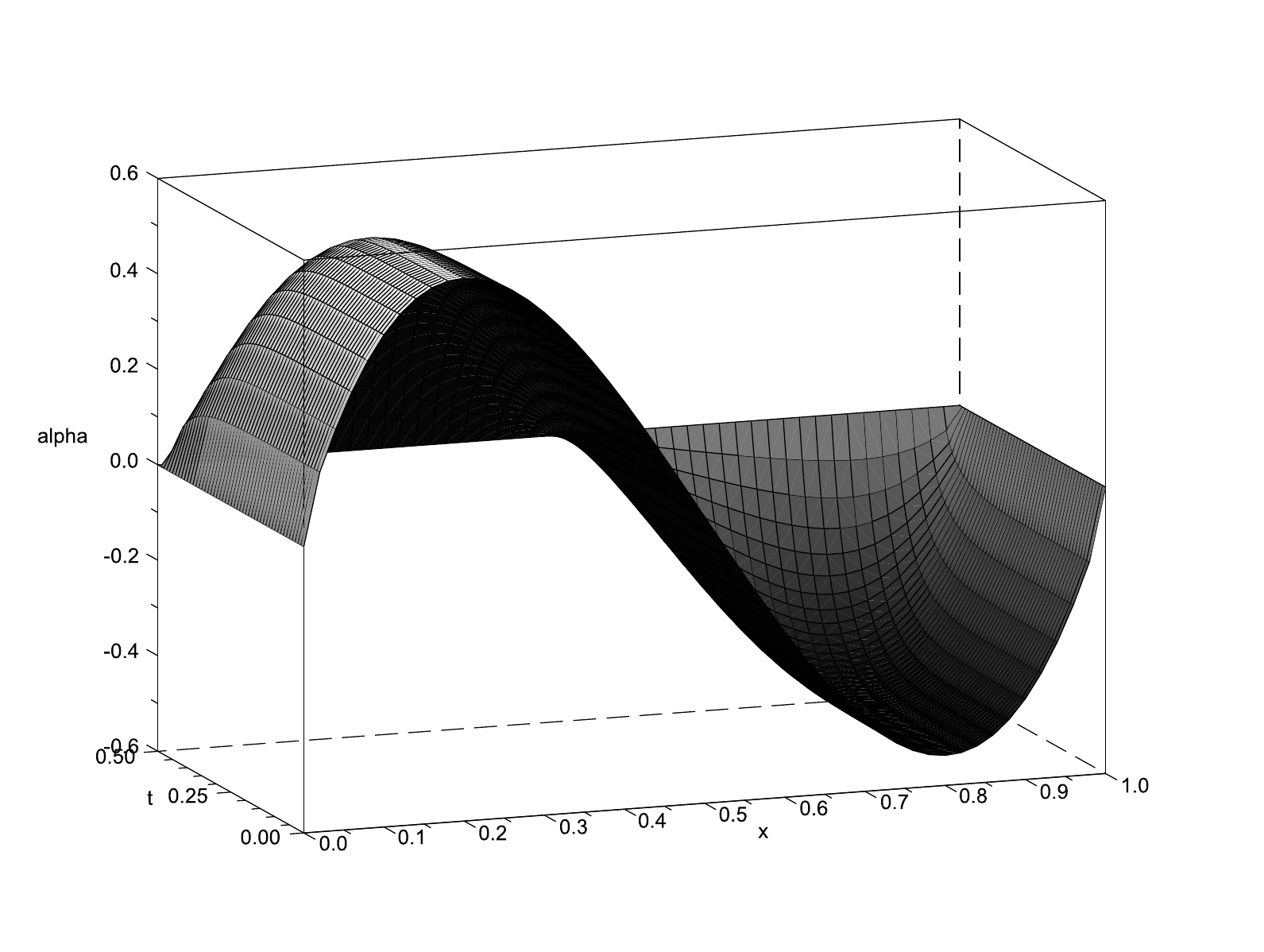}\\
  \caption{Solution for $\alpha = \nabla u$}
\end{figure}

Now, we can study, using the preceding choice of $u_T$ and $f$, the complexity of the scheme and its convergence speed.\\

We first consider the time spent by the algorithm for $\left\|\widehat{\phi}^{n+\frac 12}\widehat{\psi}^{n+1} - \widehat{\phi}^{n-\frac 12}\widehat{\psi}^{n} \right\|_{\infty}$ to be below $10^{-7}$ for the same parameters as above except $\Delta t$ and $\Delta x$ that vary:

\begin{figure}[!h]
\centering
  \includegraphics[width=260pt]{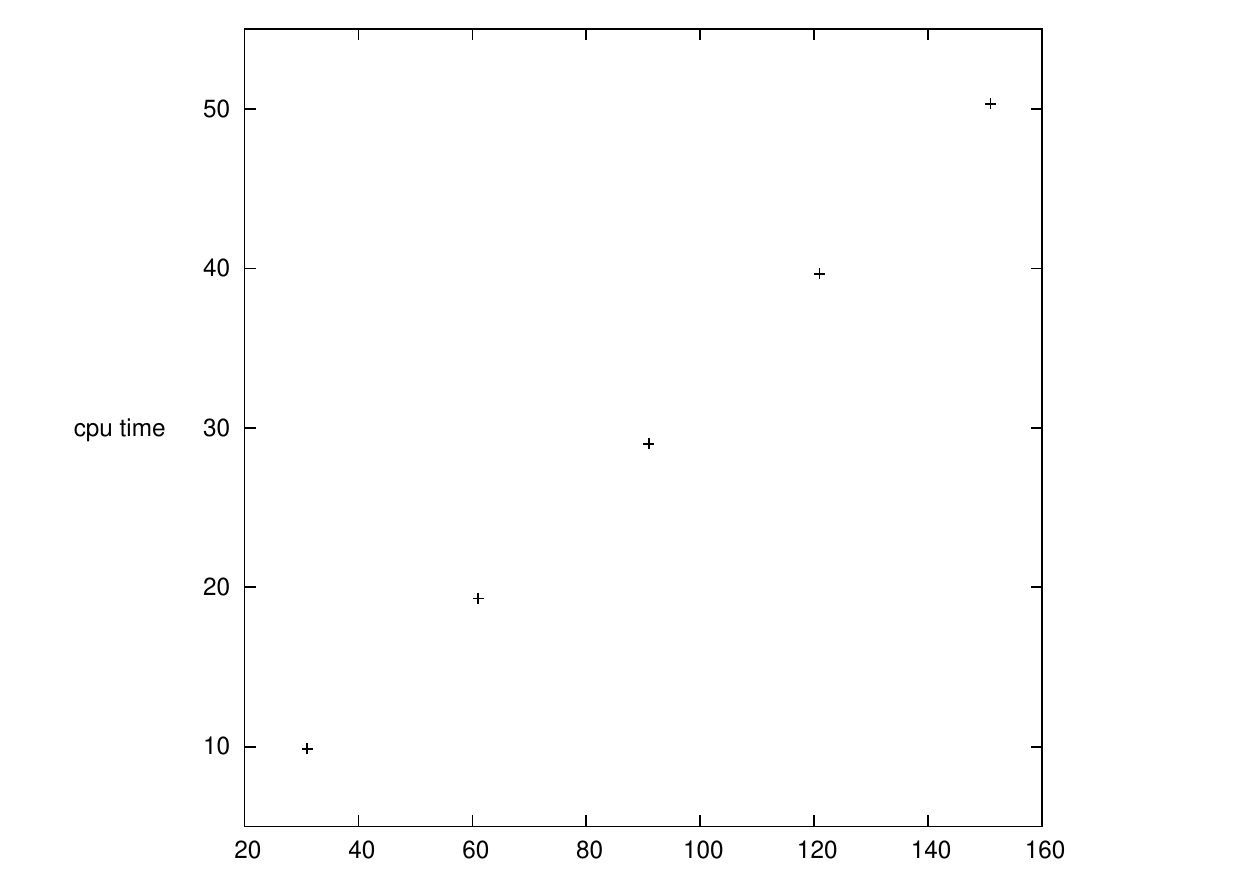}
  \hspace{0.1cm}
  \includegraphics[width=260pt]{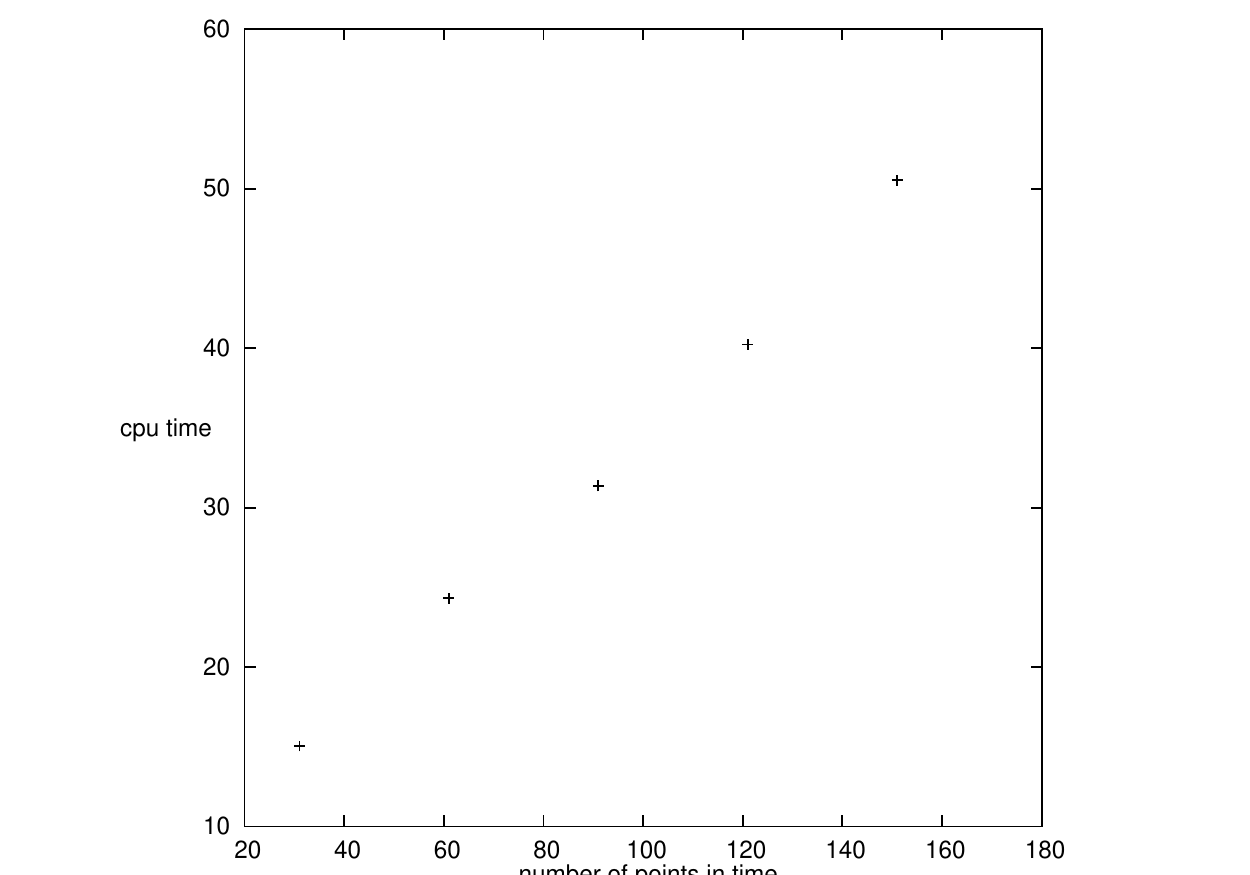}\\
  \caption{Left: Computing time for $\Delta t = 0.005$ and different values of $\Delta x$. Right: Computing time for $\Delta x = 0.01$ and different values of $\Delta t$.}
\end{figure}
\newpage

We see that the computing time is linear in the number of points in space if the number of points in time is fixed and similarly that the computing time is linear in the number of points in time if the number of points in space is fixed.\\

Now, if $\Delta t$ and $\Delta x$ are proportional, we obtain that the computing time is quadratic:

\begin{figure}[!h]
\centering
  \includegraphics[width=280pt]{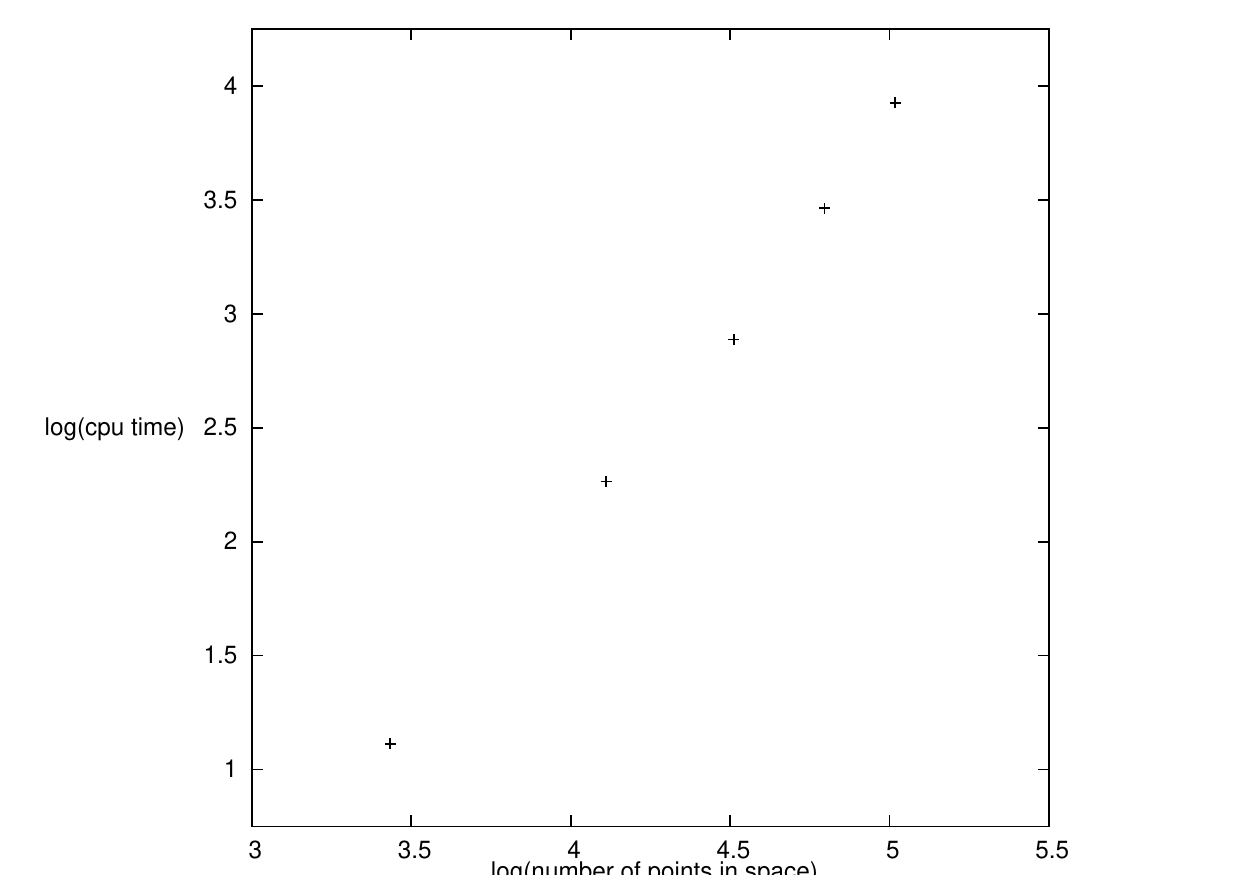}\\
  \caption{Computing time for $\Delta x = 2\Delta t$ for different values of $\Delta t$.}
\end{figure}

Empirically, we therefore have that the computing time is of order $\mathcal{O}\left((\Delta x \Delta t)^{-1}\right)$.\\

Then, we can also consider a reference solution calculated on a grid with 301 points in both space and time and find the convergence speed with respect to this reference solution. We consider first the computation of a solution -- the algorithm being stopped as soon as $\left\|\widehat{\phi}^{n+\frac 12}\widehat{\psi}^{n+1} - \widehat{\phi}^{n-\frac 12}\widehat{\psi}^{n} \right\|_{\infty}$ is below $10^{-7}$ -- for $\Delta x = \frac {1}{150}$ and different values of $\Delta t$ and calculate the error in discrete uniform norm:

\begin{figure}[!h]
\centering
  \includegraphics[width=260pt]{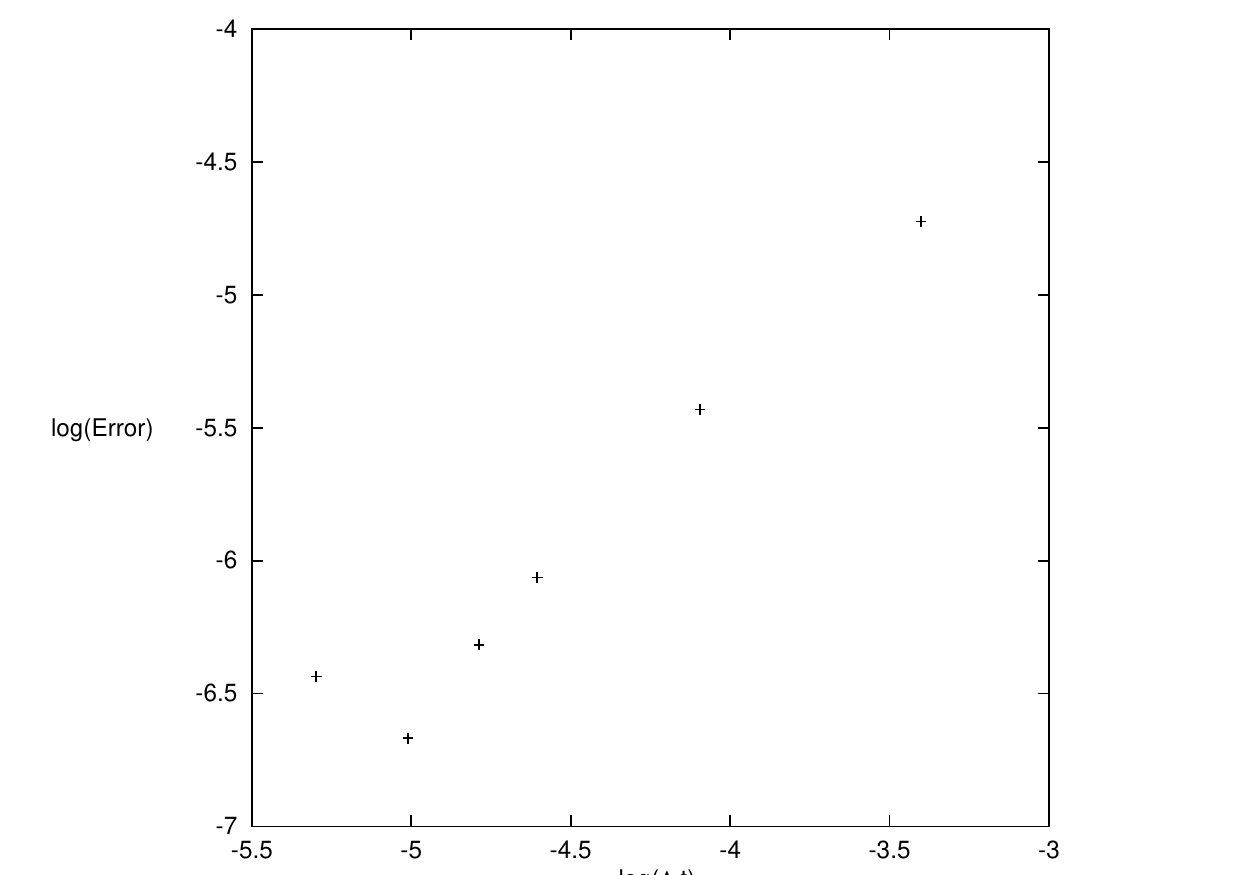}
  \hspace{0.1cm}
  \includegraphics[width=260pt]{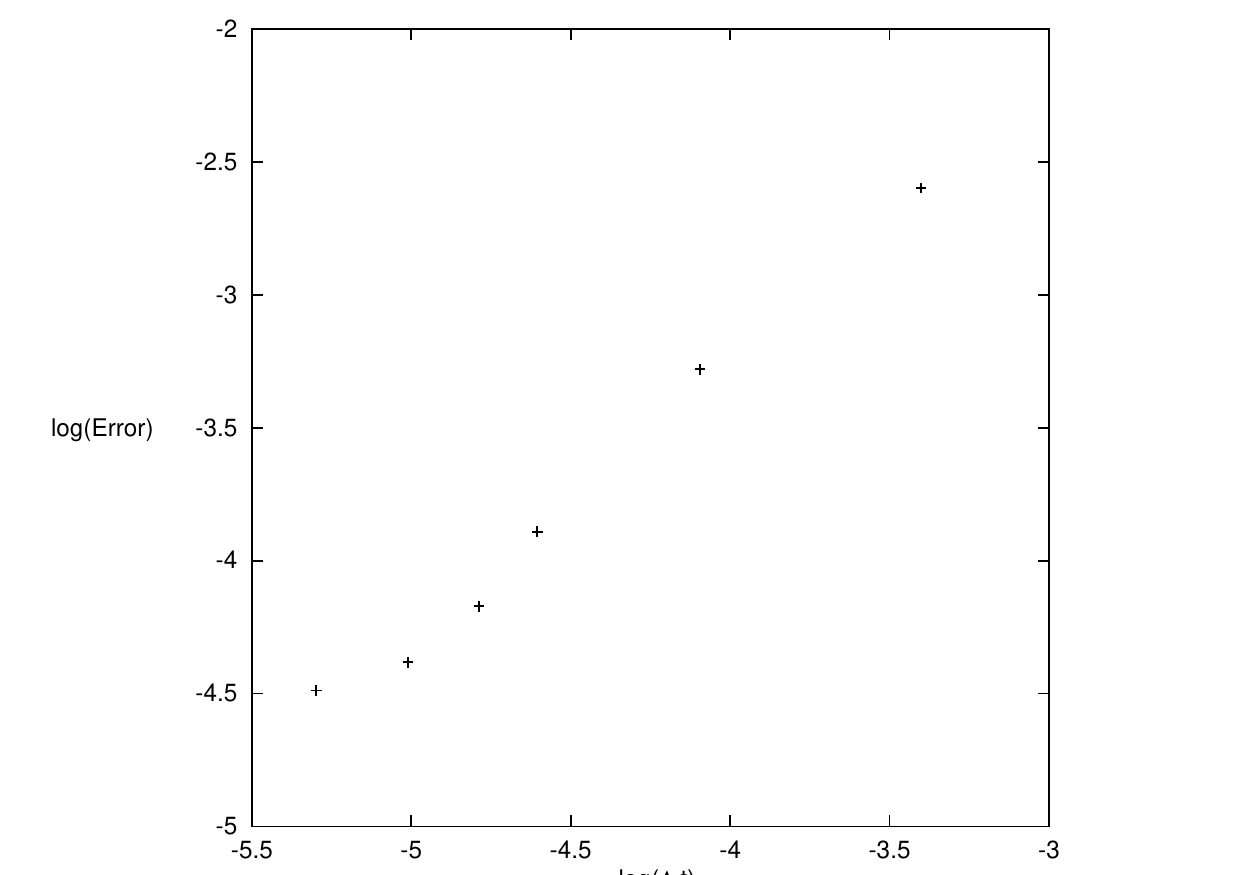}\\
  \caption{Error in discrete uniform norm for $\phi$ (left) and $\psi$ (right) with $\Delta x = \frac {1}{150}$ and different values of $\Delta t$. The reference solution has been calculated with $\Delta x = 2 \Delta t = \frac{1}{300}$.}
\end{figure}

We consider then the computation of a solution for $\Delta t = \frac {1}{300}$ and different values of $\Delta x$ and calculate the error in discrete uniform norm\footnote{The fact that the convergence speeds up for small $\Delta x$ in the second set of graphs may be due to the fact that we did not consider a closed-form reference solution but, rather, an approximation of the solution calculated on a grid.}:

\begin{figure}[!h]
\centering
  \includegraphics[width=260pt]{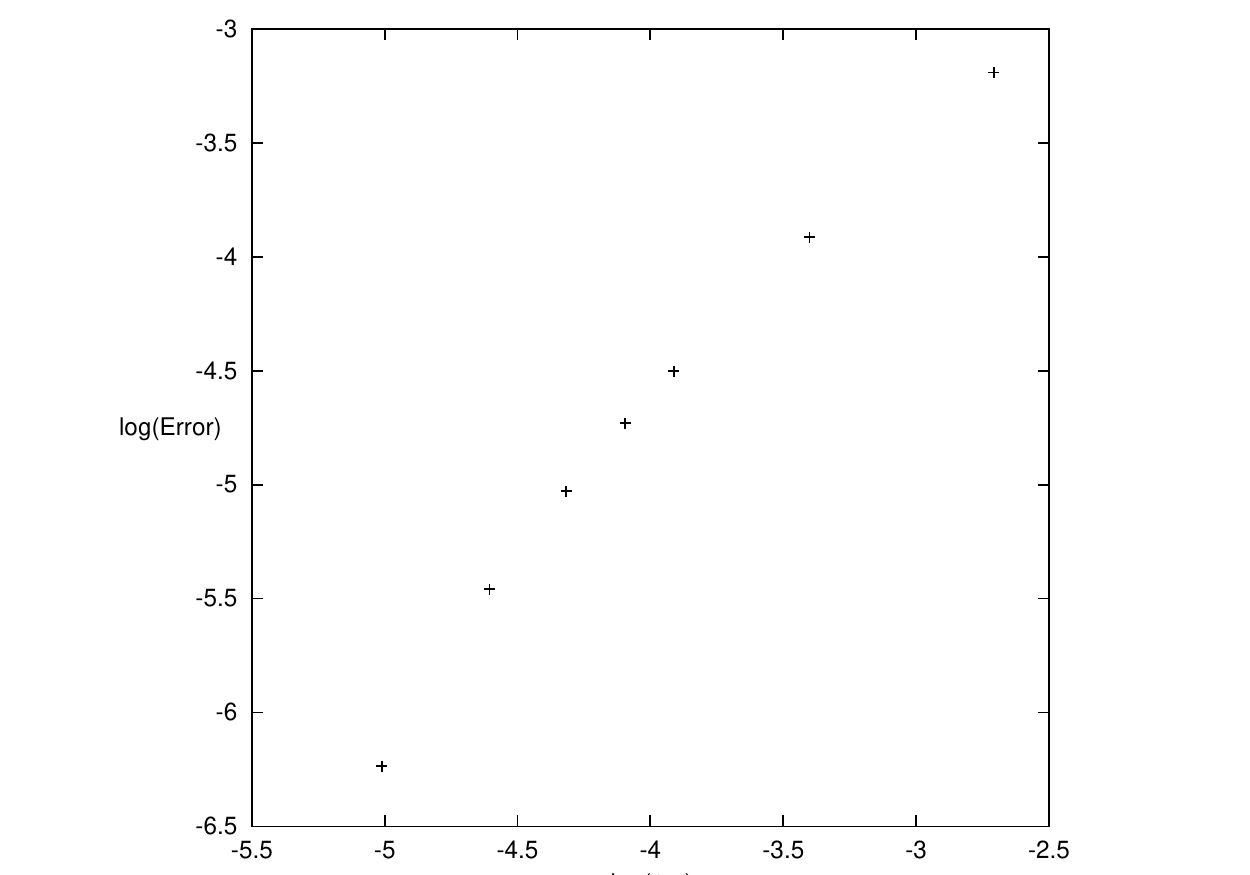}
  \hspace{0.1cm}
  \includegraphics[width=260pt]{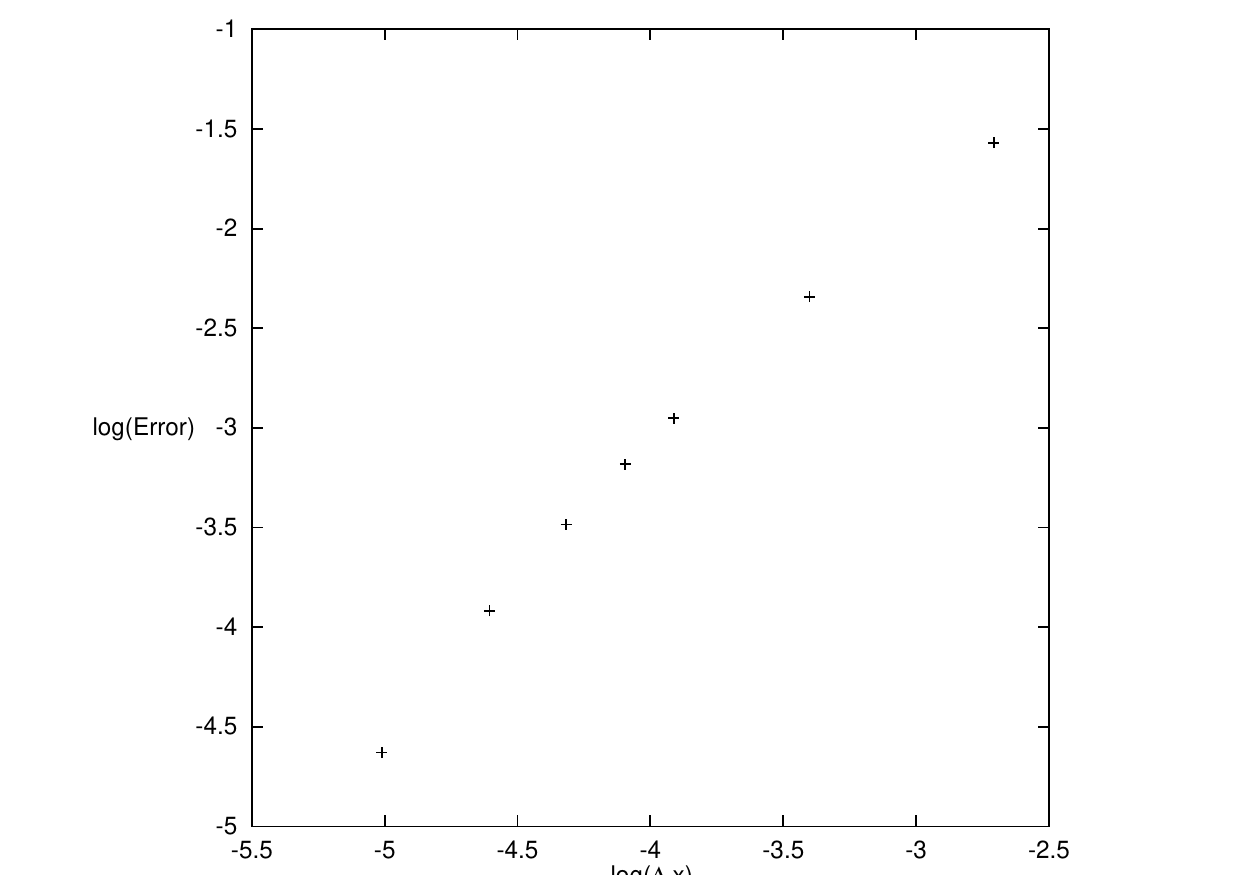}\\
  \caption{Error in discrete uniform norm for $\phi$ (left) and $\psi$ (right) with $\Delta t = \frac {1}{300}$ and different values of $\Delta x$. The reference solution has been calculated with $\Delta x = 2 \Delta t = \frac{1}{300}$.}
\end{figure}
\newpage

We see that the scheme seems to be of order $1$ in both time and space.\\

Now, we can look at the evolution of the total mass of $\widehat{m}^{n+1}$, as $n$ grows. Because of Proposition 6 and Proposition 12 we overall expect to have a loss of mass as $t$ grows, although this may not be true when $\widehat{\psi}^n$ and $\widehat{\psi}^{n+1}$ are really close, namely for large $n$. This is what we observe empirically:\\

\begin{figure}[!h]
\centering
  \includegraphics[width=240pt]{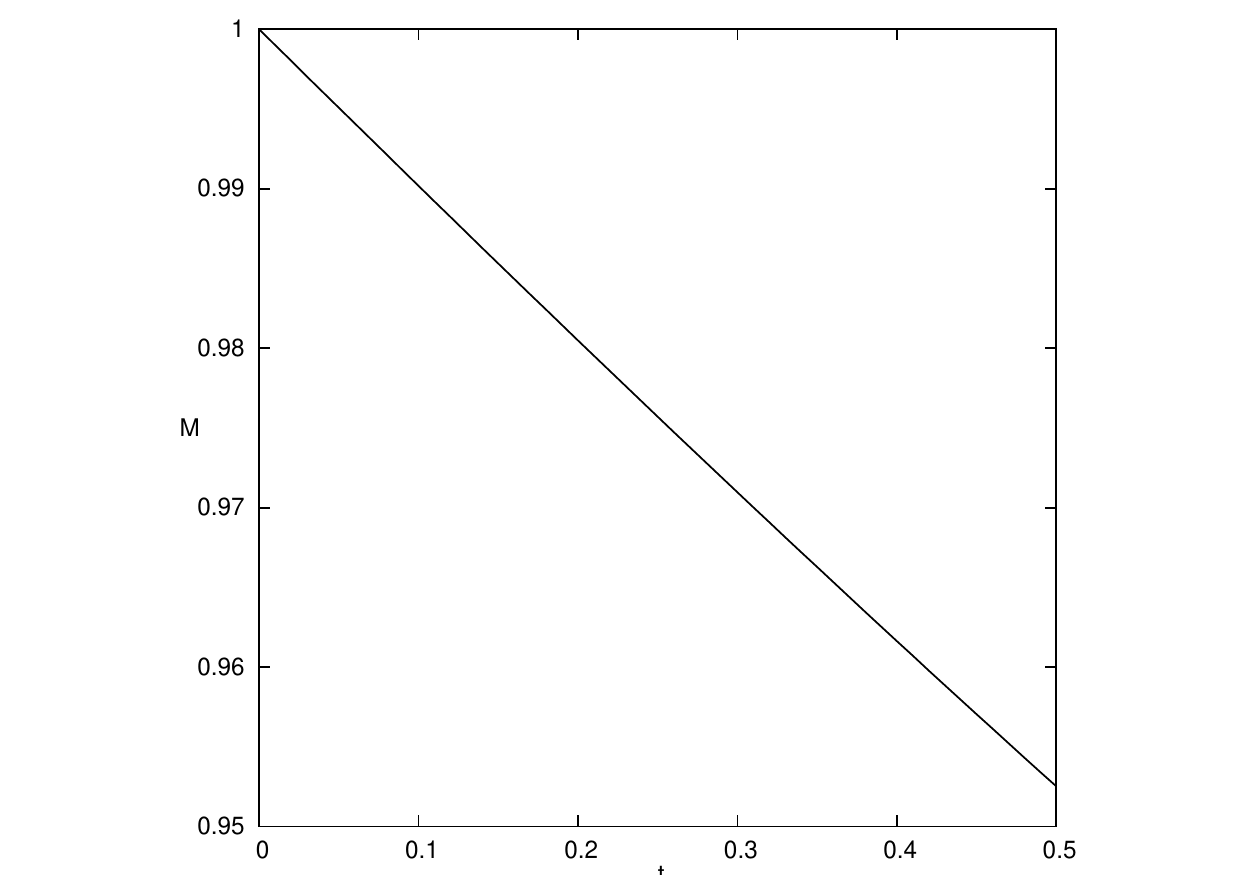}
  \hspace{0.1cm}
  \includegraphics[width=240pt]{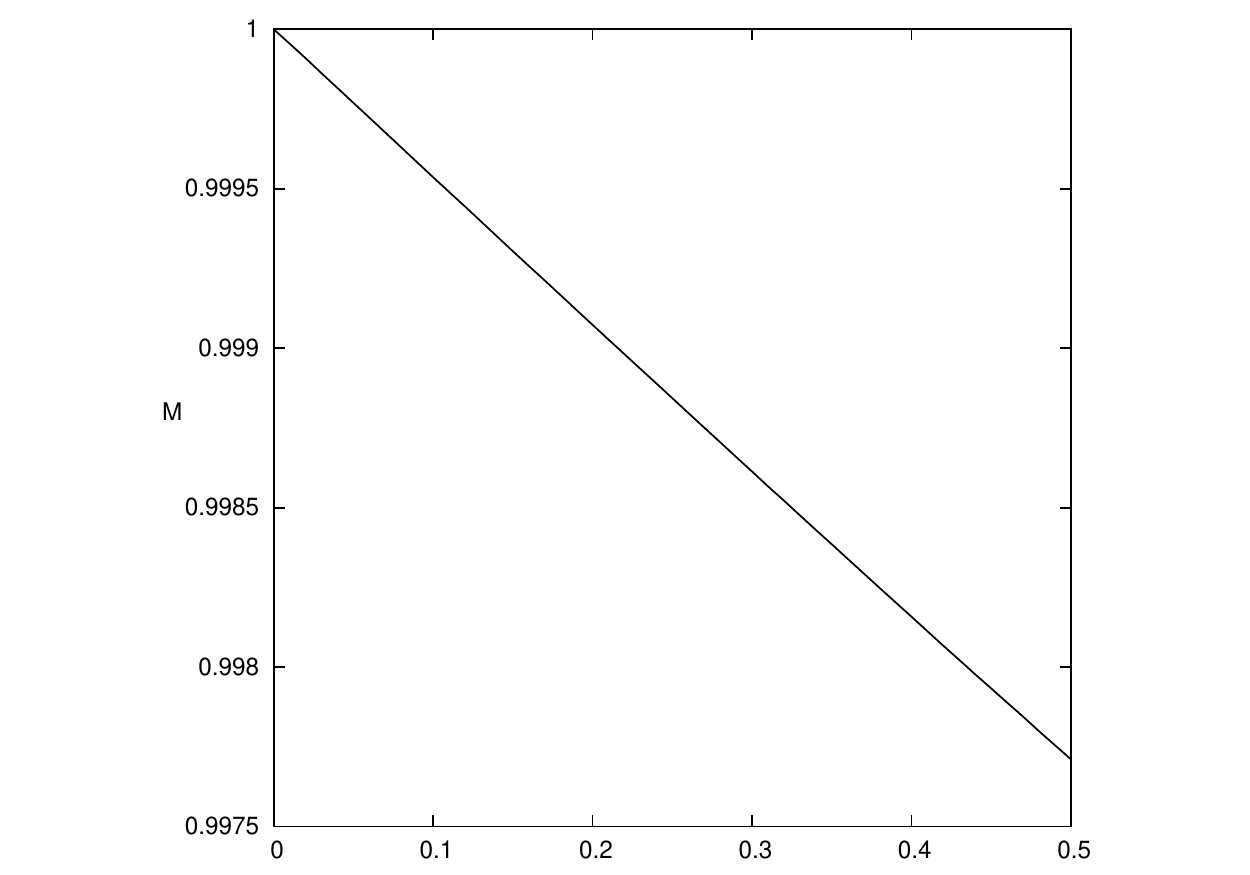}\\
  \includegraphics[width=240pt]{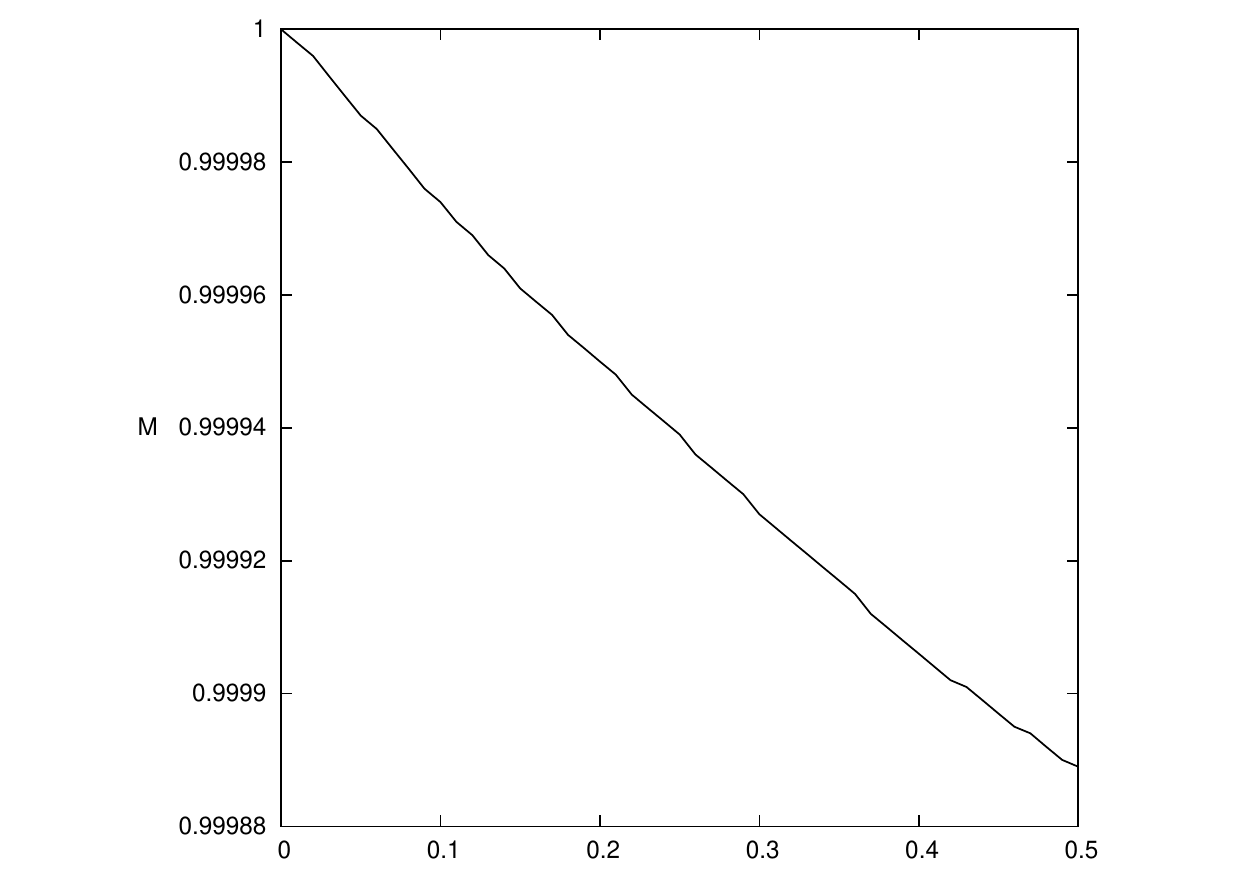}
  \hspace{0.1cm}
  \includegraphics[width=240pt]{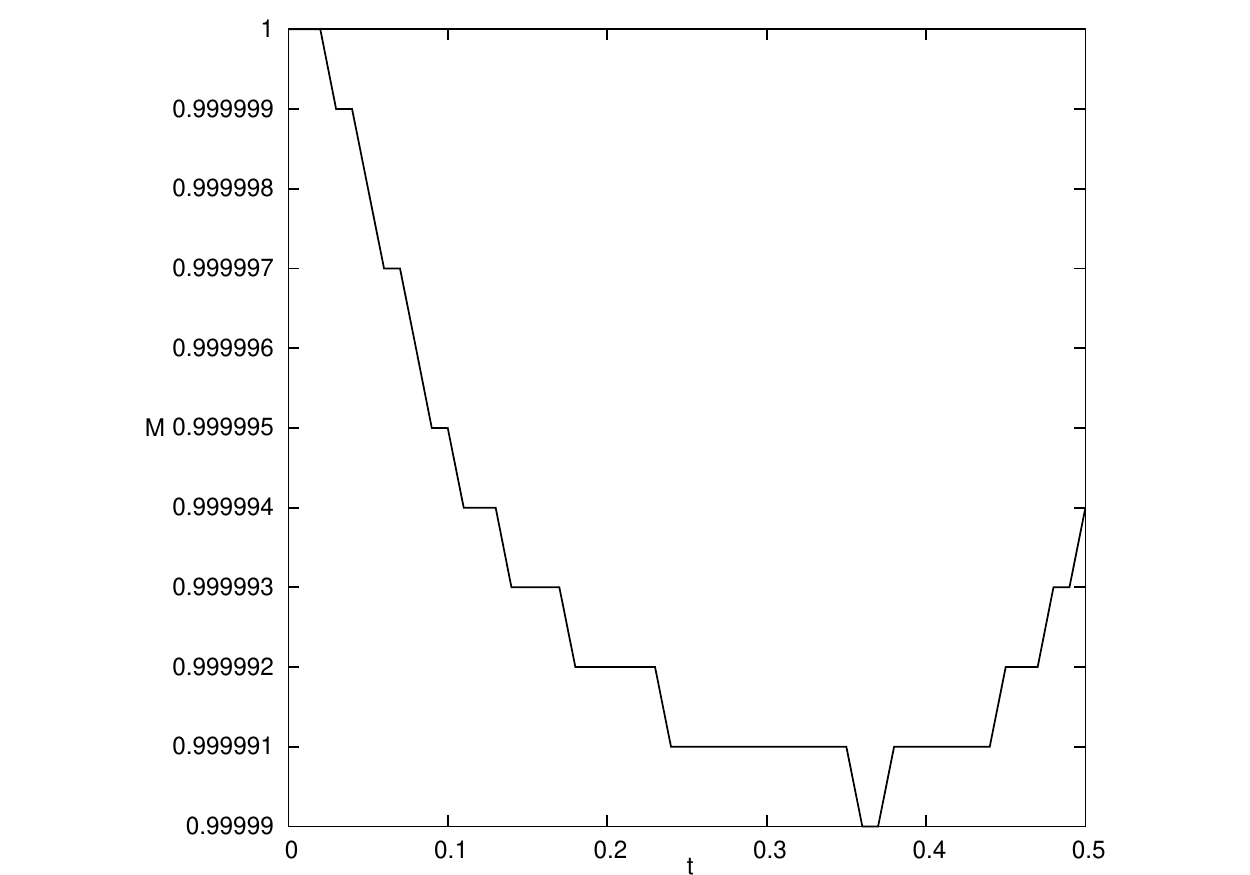}\\
  \caption{Evolution of the total mass of $\widehat{m}^{n+1}$. Top left: $n=0$, Top right; $n=1$, Bottom left: $n=2$, Bottom right: $n=3$.}
\end{figure}
\newpage
Finally, to come back to the question raised about $\sigma$, we show the computing time of the scheme for different values of $\sigma$:\\

\begin{figure}[!h]
\centering
  \includegraphics[width=320pt]{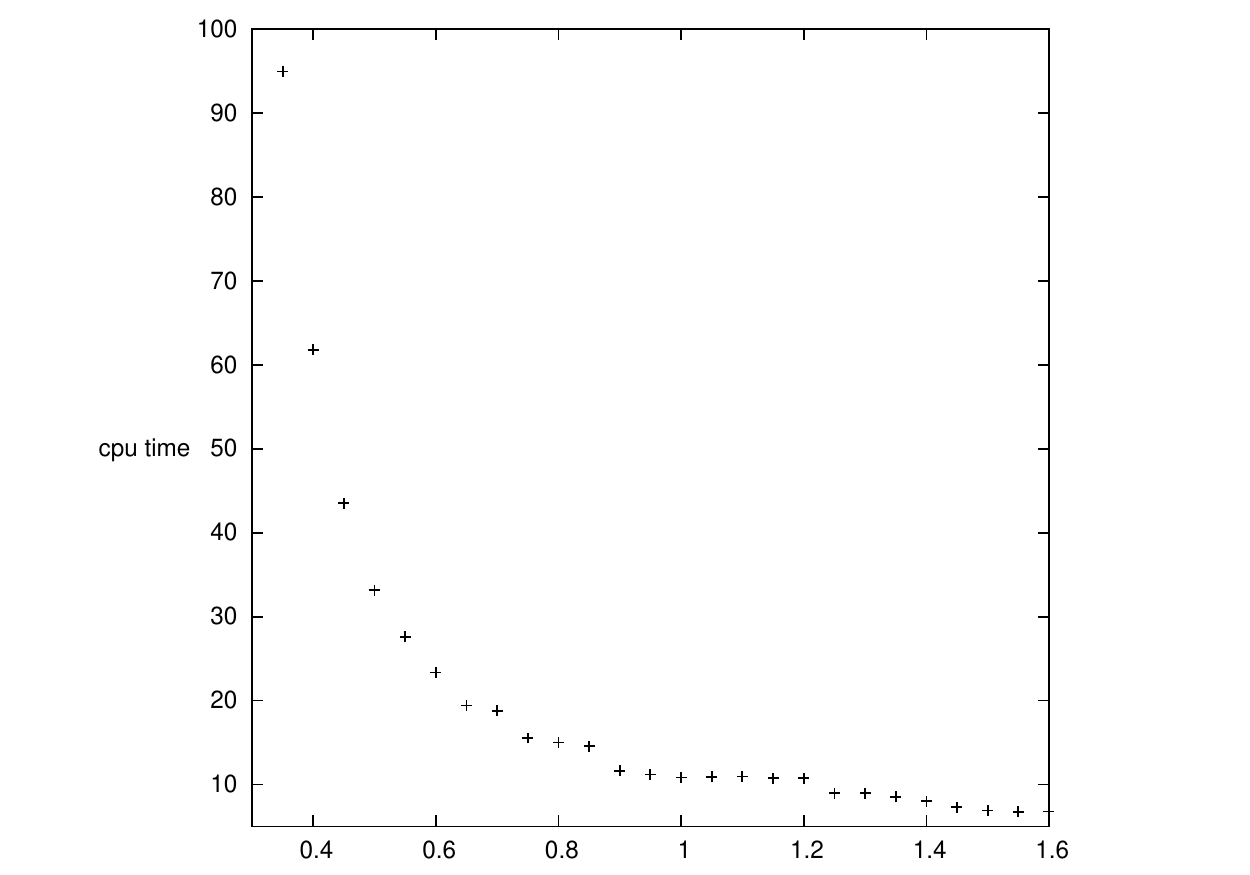}\\
  \caption{Computing time for $\Delta x = 0.02$ and $\Delta t = 0.005$ for different values of $\sigma$.}
\end{figure}

We see that the smaller the $\sigma$, the slower the algorithm. This is due to the increasing number of steps in $n$ to obtain $\left\|\widehat{\phi}^{n+\frac 12}\widehat{\psi}^{n+1} - \widehat{\phi}^{n-\frac 12}\widehat{\psi}^{n} \right\|_{\infty} < 10^{-7}$ and is also due to the increasing number of steps in the Newton scheme inside each step. The rationale behind this is that a small change in $\sigma$ does not change $u$ nor $m$ by a lot. However, $\phi = e^{\frac{u}{\sigma^2}}$ and $\psi = m e^{-\frac{u}{\sigma^2}}$ are largely influenced by $\sigma$: decreasing $\sigma$ increases $\phi$ and decreases $\psi$ by a lot.

\newpage

\bibliographystyle{plain}
\nocite{*}
\bibliography{Quadratic_Hamiltonian}

\end{document}